\documentclass[11pt]{article}
\usepackage{amsfonts}
\usepackage{oldgerm,amsmath,amssymb,theorem}
\newenvironment{proof}{{\bf Proof}:}{\vskip 5mm }

\newenvironment{rem}{{\bf Remark}:}{\vskip 5mm }

\newtheorem{lemma}{Lemma}
\newtheorem{proposition}[lemma]{Proposition}

\newtheorem{theorem}[lemma]{Theorem}

\begin{document}
\title{On the N-integrality of instanton numbers}
\author{Vadim Vologodsky}
\date{}
\maketitle
\begin{abstract}
We prove the results announced in [KSV] modulo one general fact on
Voevodsky motives that does not exist in the published literature.
Namely, we assume that the functor of motivic vanishing cycles
commutes with the Hodge and l-adic  realizations.
\end{abstract}

\begin{center}
{\large {\bf Introduction}}
\end{center}
This paper is a result of a joint work with Maxim Kontsevich and
Albert Schwarz. However, they decided not to sign it in the capacity
of authors.

 Let $\pi: X \to  C$   be a family of Calabi-Yau n-folds
over a smooth curve and let $a\in \overline C - C$ be a maximal
degeneracy point of $\pi$. We assume that the pair $( \pi: X \to C,
a)$ is defined over $\mathbb{Z}$. It is predicted that the power
series expansion for the canonical coordinate on $C$ at the point
$a$ has integral coefficients ([M]). This is a higher-dimensional
generalization of the classical fact that the Fourier coefficients
of the j-invariant are integers. This conjecture was checked in a
number of cases in [LY]. Another related conjecture says the
instanton numbers $n_d$, defined from the Picard-Fuchs equation, are
integers (see {\it loc. cit.}). We do not know how to prove these
conjectures. However, in the present paper we indicate a proof of a
weaker statement, namely, that these numbers belong to the subring
$\mathbb{Z}[N^{-1}]\subset \mathbb{Q}$ \footnote{i.e. $n_d =
\frac{m}{N^k}$, where $m$ and $k$ are integers.}, where $N$ is an
explicitly defined integer. \footnote{For example, for the quintic
family [COGP], we can take $N=2\times 3\times 5$.}

 There are two main ingredients in our proof.
The first one is {\it the Frobenius action} on the p-adic de Rham
cohomology.  It easy to see, that under our assumptions, both the
coefficients of  power series expansion for the canonical
coordinates and the instanton numbers are rational. Thus, to prove
the integrality statement, it will suffice to show that for almost
every prime $p$ they are p-adic integers. To do this we look at the
relative de Rham cohomology of our family over p-adic numbers. Then
the Frobenius symmetry or, more precisely, the existence of the
Fontaine-Laffaille structure on the cohomology bundle imply certain
strong integrality properties of the parallel sections (i.e.
solutions to the Picard-Fuchs equation ). \footnote{The idea to use
the Frobenius action is due, in a slightly different setting, to Jan
Stienstra [Sti].}

The second ingredient
 is the motivic vanishing cycles. Assume, for the purpose of Introduction, that
$dim H^n(X/C)=n+1$. Then the limit Hodge structure
 of the variation $H^n(X/C)$ is a mixed Hodge-Tate structure. Consider the
 corresponding {\it period matrix} $(a_{ij})$. This is a matrix with highly transcendental complex
 coefficients.  On the other hand, we consider the limit Fontaine-Laffaille
 module and look at the corresponding {\it Frobenius action} $(b_{ij})$. This is a
 matrix with p-adic coefficients. To complete the proof of the integrality statement we need to establish a certain relation between the superdiagonal
entries of the two
 matrices.  Namely, we have to show, that
 \begin{equation}\label{cperiod}
 a_{m\,m+1}= (2\pi i)^{m-1} log\, c
 \end{equation}
  $$
b_{m\,m+1}= \pm p^{m-1}
 log\, c^{1-p}
 $$
for some {\it rational} number $c$. \footnote{Logarithmic functions
in these formulas have different meanings: in the first formula
$log$ is the the usual complex-valued function while in the second
formula $log$ takes p-adic values.} Standard conjectures on motives
imply the existence of a mixed Artin-Tate motive $T$ over
$\mathbb{Q}$ whose Hodge and p-adic realizations are
 the limit Hodge and Fontaine-Laffaille structure correspondingly.
 This yields a certain explicit relation between {\it all} the coefficients of matrices $(a_{ij})$ and $(b_{ij})$
 and, in particular, formula (\ref{cperiod}). Unfortunately, the
 motivic conjectures needed to justify this argument are very far
 from being proved. However, we construct in Section 4 a 1-motive
 $M_{t,n}(X_{\mathbb{Q}})$  that should be thought of as the maximal
1-motive quotient of $T^*$ and then use it to prove (\ref{cperiod}).

Still, at one point we  have to rely on a general fact that has not
yet appeared in the published literature. Namely, we have to assume
the compatibility of Ayoub's motivic vanishing cycles functor with
the Hodge and l-adic realization functors. Although not published
the required compatibility is known to experts [A2], [BOV] and
hopefully this piece of a general theory will be written in a matter
of time.

The paper is organized as follows. Section 1 contains statements of
the results. In Section 2 we recall some well known facts on vector
bundles with logarithmic connection, variations of Hodge structure,
and give an interpretation of the canonical coordinate as an
extension class of certain variations of Hodge structure. In Section
3 we define, using p-adic Hodge Theory,
 a p-adic analog of the canonical coordinate and Yukawa coupling (for 1-parameter families of Calabi-Yau varieties over $\mathbb{Z}_p$)
and prove, using Dwork's Lemma, the (p-adic) integrality statements
for these objects. Finally, in the last section (the most technical
one) we show for families over $\mathbb{Z}$ that the two
constructions (complex and p-adic) give the same functions.
\footnote{This amounts to proving relation (\ref{cperiod}).} To do
this we give a third geometric definition (i.e. which makes sense
over any ground field) of the canonical coordinate.
 The construction is based on the notion of motivic nearby cycles (due to Ayoub [A1])
 and uses the language of Voevodsky motives.

 {\bf Acknowledgments.} Besides Maxim Kontsevich and Albert
Schwarz whom many ideas in the present paper are due to, I wish to
express my gratitude to Joseph Ayoub, who drew my attention to the
paper [BK]. Special thanks go to the referee for his generous help
in turning a raw draft into a paper. In particular, the referee
pointed out that the compatibility of the motivic Albanese functor
with the Hodge realization, that plays an important role in our
argument, is not proved in the published literature. At the end the
author wrote a separate paper ([Vol]) with a proof of the required
compatibility. \\
I acknowledge the stimulating atmosphere of IHES and MPIM in Bonn
where parts of this work were done.\\
This research was partially supported by NSF grant DMS-0401164.

  \section{Terminology and statements of the results}

{\bf 1.1.  Definition of the canonical coordinate. }  The following
construction is due to Morrison [M].
 Let $\pi: X \to  C$   be a family of Calabi-Yau varieties of dimension $n$
over a smooth curve over $\mathbb{C}$. This means that locally on
$C$ the relative canonical bundle $\Omega^n_{X/C}$ is trivial.
Assume that $C$ is embedded into a larger smooth curve  $\overline C
\supset C $ and $a \in \overline C - C$  is a boundary point. The
point $a$ is called a maximal degeneracy point if the monodromy
operator $M:  H_n(X_{a'}, \mathbb{Q}) \to H_n(X_{a'}, \mathbb{Q})$,
corresponding to a small loop around $a$ is unipotent and $(M -
Id)^n \ne 0$.

{\bf Remark.}
 It is known that for any smooth proper family $\pi: X \to  C$ with a unipotent monodromy,  $(M - Id)^{n+1} = 0$.

 From now on we will assume that $a$ is a  maximal degeneracy point.
 Set $N_B= log\, M$. The following remarkable result was derived by Morrison from the very basic properties of the
 limit Hodge structure \footnote{The proof
is reproduced in 2.2.}.

\begin{lemma}([M], Lemma 1.) \label{M} $dim_{\mathbb{Q}} Im\, N_B^n =1$ and $dim_{\mathbb{Q}} Im\, N_B^{n-1} =2$.
\end{lemma}

 Denote by  $ {\cal T}_{\mathbb{Z}}$ the local system over a punctured neighborhood $D^*$ of $a$, whose fiber over a point $a'\in C$
is $Im( H_n(X_{a'}, \mathbb{Z}) \to  H_n(X_{a'}, \mathbb{Q})) $. Let
$\delta_1, \delta_2$ be a basis  for $Im\, N_B^{n-1} \cap  ({\cal
T}_{\mathbb{Z}})_{a'} $ such that
  $\delta_1 $  generates    $Im\, N_B^{n} \cap  ({\cal T}_{\mathbb{Z}})_{a'} $ and such that $N_B(\delta_2)$ is a positive multiple of
$\delta_1$. We may view $\delta_1$ as a section of ${\cal
T}_{\mathbb{Z}} $ over $D^*$ and $\delta_2$ as a section of  the
quotient ${\cal T}_{\mathbb{Z}}/\mathbb{Z}\delta _1 $. Choose a
non-vanishing section $\omega $ of  $\pi_*\Omega^n_{X/C}$ over
$D^*$.
 We then see that
\begin{equation}\label{cc}
  q= exp(2\pi i \frac{\int _{\delta _2} \omega}{ \int _{\delta _1} \omega})
\end{equation}
is a well defined function on a punctured neighborhood $ D^*$ and
that $q$ does not depend on the choice of $\delta_i$ and $\omega$ we
made. Moreover, the function $q$ extends to $a$ and $ord_a  q = k$,
where $k$ is defined from the equation $N_B(\delta_2)=  k \delta
_1$. \footnote{ The nontrivial part is to show that ${ \int _{\delta
_1} \omega}$ does not vanish on a sufficiently small $D^*$ and that
that $q$ has regular singularity at $a$. This is another corollary
of the existence of the limit Hodge structure.
 See 2.2.}
We shall say  the Betti monodromy of the family $X\to C$ is small if
$k=1$. In this case, the function $q$ is called  the canonical
(local) coordinate on $\overline C$.

{\bf 1.2. Yukawa function.} Denote by
$${\cal H}= ({\cal H}_{\mathbb{Z}}= Im( R^n \pi_* \mathbb{Z} \to  R^n \pi_* \mathbb{Q}),\,  {\cal F}^{n} \subset {\cal F}^{n-1}
\subset \cdots \subset {\cal F}^0 = {\cal H}_{\mathbb{Z}} \otimes
{\cal O}_{D^*})$$
 the variation of Hodge structure associated to $\pi: X\to C$.
The  Kodaira-Spencer operator
$$\overline \nabla: {\cal F}^p/{\cal F}^{p+1} \to {\cal F}^{p-1}/{\cal F}^p \otimes \Omega^1_{D^*} $$
 extends to a homomorphism of graded algebras
  $$ S^{\cdot}T_{D^*} \to End_{{\cal O}_{D^*}}(\bigoplus_ p {\cal F}^p/{\cal F}^{p+1}).$$
Specializing, we get a morphism
\begin{equation}\label{ks}
\kappa :  S^{n}T_{D^*} \to Hom _{{\cal O}_{D^*}} ({\cal F}^n, {\cal
F}^0/{\cal F}^1)\simeq ({\cal F}^n \otimes {\cal F}^n)^*.
\end{equation}
The line bundle ${\cal F}^n \otimes {\cal F}^n$ is naturally
trivialized over $D^*$. To see this, denote by
 $\omega \in {\cal F}^n$  the differential form such that
$$\int _{\delta _1} \omega = (2\pi i)^n. $$
The section $\omega \otimes \omega \in  {\cal F}^n \otimes {\cal
F}^n$ defines the desired trivialization. \footnote{The cycle
$\delta_1$ is defined up to sign. But the trivialization of ${\cal
F}^n\otimes {\cal F}^n$ is independent of this choice.} Define the
Yukawa function on $D^*$ to be
$$Y = \kappa((q\frac{d}{dq})^n)\cdot (\omega \otimes \omega).$$
One can check that $Y$ extends to $D$.\footnote{This is a corollary
of a result of Schmid, which says that the Hodge filtration extends
to Deligne's canonical extension of the underlying vector bundle.
(See also 2.1).}

{\bf 1.3. Statement of main results.} Let $S= spec \,
\mathbb{Z}[N^{-1}]$ be an open subscheme of $spec \, \mathbb{Z}, \;
$ $\overline C_S$ a smooth curve over $S$, and let $a: S
\hookrightarrow \overline C_S$ be a section.  Denote by $t$
 a local coordinate on a open neighborhood of $a$ such that $t(a)=0$.
 Let $\pi: X_S \to C_S = \overline C_S- a$ be a smooth proper family of Calabi-Yau
schemes. We will make the following assumptions: \\
i) $a_{\mathbb{C}} $ is the maximal degeneracy point of the complex
family $X_{\mathbb{C}}\to C_{\mathbb{C}}$ \\
ii) $\pi:  X_S \to C_S $  extends to a semi-stable morphism
 $\overline{\pi}: \overline X_S \to \overline C_S$ \footnote{ i.e. locally for the etale topology
$\overline{\pi}: \overline X_S \to \overline C_S$ is isomorphic to
$spec\, \mathbb{Z}[t, x_1, \cdots x_n]/(x_1\cdots x_r - t) \to
spec \, \mathbb{Z}[t]$, where $r\leq n $ .}\\
iii) All primes $p\leq dim \, X_{\mathbb{C}}$ are invertible on $\mathbb{Z}[N^{-1}]$ \\
iv)  The Betti monodromy of the family  $X_{\mathbb{C}}\to
C_{\mathbb{C}}$ is small ( see  1.1 ).

\begin{theorem}\label{th1} Assume that $q'(0)$ is a rational number . Then
 $$q(t)\in  (\mathbb{Z}[N^{-1}]((t)))^* . $$
 \end{theorem}
\begin{rem} We shall see in Section {\bf 4.5} that, for any
family $X_{\mathbb{Q}} \to C_{\mathbb{Q}} $ over $\mathbb{Q}$ with a
maximal degeneracy point at $a\in \overline
C_{\mathbb{Q}}(\mathbb{Q})$, $q'(t)^r\in \mathbb{Q}((t))$, for some
integer $r$.
\end{rem}

For the next result we assume that $dim \, X_{\mathbb{C}} = 4$ (i.e.
$X_{\mathbb{C}}\to C_{\mathbb{C}}$ is a family of threefolds) and
\footnote{Observe that dimension of the space of first order
deformations of a Calabi-Yau n-fold Y is equal to $dim \, H^1(Y,
T_Y)= dim \, H^{n-1}(Y, \Omega ^1)$. Thus the  condition (\ref{one})
implies that $\pi: X_{\mathbb{C}}\to C_{\mathbb{C}}$ induces a
dominant map from $C_{\mathbb{C}}$ to an irreducible component of
the moduli space of Calabi-Yau threefolds.
 The case of a higher dimensional component in the moduli space of Calabi-Yau threefolds will be considered elsewhere (also see [KSV], Section 3). }
\begin{equation}\label{one}
rk \, H^3_{DR}(X_{\mathbb{C}}/C_{\mathbb{C}}) = 4.
\end{equation}
We also assume that $q'(0)$ is a rational number.
\begin{theorem}\label{th2} One has
\begin{equation}\label{instanton}
Y(q)= n_0+ \sum^{\infty}_{d=1}n_d d^3 \frac{q^d}{1-q^d},
\end{equation}
where $n_i \in \mathbb{Z}[N^{-1}]$. \footnote{One readily sees that
any power series $Y(q)\in \mathbb{C}[[q]]$ can be written in the
form
 (\ref{instanton}), with $n_d\in \mathbb{C}$. Thus the content of the theorem is the integrality property of the numbers.
These are the instanton numbers the title of our paper refers to.}
\end{theorem}

\section{Hodge Theory}

{\bf 2.1. Logarithmic De Rham cohomology.} We will need the
following construction from [Ste]. Let $\overline{\pi}: \overline
X_S \to \overline C_S$ be a semi-stable morphism.

   $$
\def\normalbaselines{\baselineskip20pt
\lineskip3pt  \lineskiplimit3pt}
\def\mapright#1{\smash{
\mathop{\to}\limits^{#1}}}
\def\mapdown#1{\Big\downarrow\rlap
{$\vcenter{\hbox{$\scriptstyle#1$}}$}}
\begin{matrix}
 Y_S= \overline X_S \times _{\overline C_S} S   &    \stackrel{}{\hookrightarrow}
 &  \overline X_S & \hookleftarrow  & X_S \cr
\mapdown{}  &  &\mapdown{\overline{\pi}} & & \mapdown{\pi}    \cr
 S &    \stackrel{a}{\hookrightarrow}
 &  \overline C_S & \hookleftarrow & C_S
 \end{matrix}
 $$
We then consider the relative logarithmic De Rham complex
$(\Omega^*_{\overline X_S/\overline C_S}(log \,Y_S), d)$ on
$\overline X_S$ defined as follows. Let  $\Omega^i_{\overline
X_S}(log \,Y_S)$ be the sheaf of differential forms with logarithmic
singularities along $Y_S$, and let $\Omega^*_{\overline
X_S/\overline C_S}(log \,Y_S)$ be the quotient of the sheaf of
algebras   $\Omega^*_{\overline X_S}(log \,Y_S)$ by the ideal
generated by $\overline {\pi}^* \eta, \, \, \eta \in
\Omega^1_{\overline C_S}(log \, S)$. One immediately sees that
$\Omega^*_{\overline X_S/\overline C_S}(log \,Y_S)$ is a locally
free sheaf of ${\cal O}_{\overline X_S}$-modules and that the
exterior differential $d: \Omega^i_{\overline X_S}(log \,Y_S)\to
 \Omega^{i+1}_{\overline X_S}(log \,Y_S)$ descends to $\Omega^*_{\overline X_S/\overline C_S}(log \,Y_S)$ . We then define the logarithmic De Rham
cohomology by
  $$H^i_{log}(\overline X_S/\overline C_S)= R^i\overline{\pi}_*(\Omega^*_{\overline X_S/\overline C_S}(log \,Y_S), d).$$
$H^i_{log}(\overline X_S/\overline C_S)$ is a coherent sheaf on
$\overline C_S$ equipped with a logarithmic connection:
$$\nabla: H^i_{log}(\overline X_S/\overline C_S)\to H^i_{log}(\overline X_S/\overline C_S) \otimes \Omega^1_{\overline C_S}(log \, S).$$
Assume that all primes $p\leq n= dim \, _{C_S} X_S$ are invertible in $\mathbb{Z}[N^{-1}]$. Then \\
i) the coherent sheaf $H^i_{log}(\overline X_S/\overline C_S)$ is
locally isomorphic to a direct sum of the sheaves  ${\cal
O}_{\overline C_S}$,
${\cal O}_{\overline C_S}/p^e$. In particular, the quotient of $H^i_{log}(\overline X_S/\overline C_S)$ modulo torsion is locally free.\\
ii) the Hodge spectral sequence (i.e. the spectral sequence
associated to the stupid filtration on $(\Omega^*_{\overline
X_S/\overline C_S}(log \,Y_S), d)$) degenerates in the first term.
Moreover the induced filtration $\overline{{\cal F}}^{\cdot}\subset
H^i_{log}(\overline X_S/\overline C_S)$ splits (in the category of
${\cal O}_{\overline C_S}$-modules) locally on $\overline C_S$. In
particular, $R^0\overline{\pi}_*(\Omega^i_{\overline X_S/\overline
C_S}(log \,Y_S))$ is a locally free
${\cal O}_{ \overline C_S}$-module. \\
iii) the residue of the connection
  $$N_{DR}= Res \, \nabla: a^* H^i_{log}(\overline X_S/\overline C_S) \to  a^* H^i_{log}(\overline X_S/\overline C_S) $$
is nilpotent. In particular, $H^i_{log}(\overline
X_{\mathbb{Q}}/\overline C_{\mathbb{Q}})$ is the Deligne canonical
extension of the vector bundle
$ H^i_{DR}(X_{\mathbb{Q}}/ C_{\mathbb{Q}})$ equipped with the Gauss-Manin connection. \\
iv) there is a canonical pairing
\begin{equation}\label{derhampairing}
<\cdot ,\cdot >_{DR}:  H^i_{log}(\overline X_S/\overline C_S)\otimes
H^{2n -i}_{log}(\overline X_S/\overline C_S) \to
 H^{2n}_{log}(\overline X_S/\overline C_S)\simeq {\cal O}_{\overline C_S}
\end{equation}
The induced pairing on the quotient of $H^*_{log}(\overline X_S/\overline C_S)$ modulo torsion is perfect.\\
Over $\mathbb{C}$ these facts are proven in [Ste]; the integral
version is contained in [Fa] (Theorems 2.1 and 6.2).

The De Rham isomorphism
$$ H^i_{DR}(X_{\mathbb{C}}/ C_{\mathbb{C}})  \simeq
R^i \pi^{an}_{\mathbb{C}*} \mathbb{Z} \otimes_{\mathbb{Z}}    {\cal
O}_{C_{\mathbb{C}}}$$ and a choice of a local coordinate $t$ on
$\overline C_{\mathbb{C}}$ yield  an integral structure (of a
topological nature) on the vector space $a_{\mathbb{C}}^*
H^i_{log}(\overline X_{\mathbb{C}}/\overline C_{\mathbb{C}})$:
\begin{equation}\label{derhamv}
a_{\mathbb{C}}^* H^i_{log}(\overline X_{\mathbb{C}}/\overline
C_{\mathbb{C}})\simeq \Psi_t^{an, un}(R^i \pi^{an}_{\mathbb{C}*}
\mathbb{Z}) \otimes _{\mathbb{Z}}    \mathbb{C}
\end{equation}
To see this, let  $Log^{\infty}= {\cal O}_{\mathbb{C}^*}[log\, t]$
be the universal unipotent local system on $\mathbb{C}^*$ and let
$\overline{Log}^{\infty}={\cal O}_{\mathbb{C}}[log\, t] $ be the
Deligne extension of  $Log^{\infty}$ to $\mathbb{C}$.  Let us view
$Log^{\infty}$ as a subsheaf of the direct image $(exp)_* {\cal
O_{\mathbb{C}}}$ of the structure sheaf on the universal cover $exp:
\mathbb{C} \to \mathbb{C}^*$. Define a $\mathbb{Z}$-lattice
$Log^{\infty}_{\mathbb{Z}} \subset Log^{\infty}$  to be $(exp)_*
\mathbb{Z} \cap Log^{\infty}$.  We have then
$$Log^{\infty}_{\mathbb{Z}}\otimes {\cal O}_{\mathbb{C}^*} \simeq Log^{\infty}.$$
 Let  $a_{\mathbb{C}} \in D \subset \overline C_{\mathbb{C}}(\mathbb{C})$ be
a disk such that the map $t: D \hookrightarrow \mathbb{C}$ defined
by the coordinate is an embedding. Given a $\mathbb{Z}$-local system
${\cal H}_{\mathbb{Z}}$ over $D^*$ we define the space of unipotent
vanishing cycles by
$$\Psi_t^{an, un}({\cal H}_{\mathbb{Z}}):= H^0(D^*,  {\cal H}_{\mathbb{Z}}\otimes
  (t|_{D^*})^* Log^{\infty}_{\mathbb{Z}} ).
\footnote{This definition is borrowed from [B].}$$ Assume that
${\cal H}_{\mathbb{Z}}$ is unipotent and denote by $\overline {\cal
H}$ the Deligne extension of  ${\cal H}={\cal H}_{\mathbb{Z}}\otimes
{\cal O}_{D^*}$ to $D$. We shall define a canonical isomorphism:
\begin{equation}\label{BS}
a_{\mathbb{C}}^* \overline {\cal H}   \simeq \Psi_t^{an, un}( {\cal
H}_{\mathbb{Z}}) \otimes _{\mathbb{Z}}    \mathbb{C}.
\end{equation}
This will induce (\ref{derhamv}). To construct (\ref{BS}) observe
that for any vector bundle $E$  over $D$
 with a logarithmic nilpotent connection (i.e. a logarithmic connection such that $N_{DR}$ is nilpotent)
we have
\begin{equation}\label{ncon}
(E|_{D^*})^{\nabla}\stackrel{\simeq}{\longleftarrow}
E^{\nabla}\stackrel{\simeq}{\longrightarrow}
    (a_{\mathbb{C}}^* E)^{N_{DR}=0}.
\end{equation}
We apply (\ref{ncon}) to the ind-object $\overline {\cal H} \otimes
 t^* (\overline{Log}^{\infty})$ and use a canonical isomorphism
$$  (a_{\mathbb{C}}^* \overline {\cal H} \otimes a_{\mathbb{C}}^*  t^* (\overline{Log}^{\infty}))^{N_{DR}=0}
\simeq a_{\mathbb{C}}^* \overline {\cal H} ,$$ which takes an
element $ v\otimes log^k t \in  (a_{\mathbb{C}}^* \overline {\cal H}
\otimes a_{\mathbb{C}}^*  t^* (\overline{Log}^{\infty}))^{N_{DR}=0}
$ to $v$ if $k=0$ and to $0$ otherwise.

 Denote by $$N_B:
\Psi_t^{an, un}(R^i \pi^{an}_{\mathbb{C}*} \mathbb{Z}) \otimes
_{\mathbb{Z}} \mathbb{Q} \to \Psi_t^{an, un}(R^i
\pi^{an}_{\mathbb{C}*} \mathbb{Z}) \otimes _{\mathbb{Z}} \mathbb{Q}
$$ the logarithm of the monodromy operator and by
$$< \cdot, \cdot>_B: \Psi_t^{an, un}(R^i \pi^{an}_{\mathbb{C}*} \mathbb{Z})\otimes
\Psi_t^{an, un}(R^{2n-i} \pi^{an}_{\mathbb{C}*} \mathbb{Z}) \to
\mathbb{Z} $$ the pairing induced by the Poincare duality. We then
have
 $$N_B= -2\pi i N_{DR} \, , \, \,  < \cdot, \cdot>_B= (2\pi i)^n < \cdot, \cdot>_{DR}.  $$

{\bf 2.2. A variation of mixed Hodge structure.} Let
$\pi_{\mathbb{C}}: X_{\mathbb{C}} \to C_{\mathbb{C}}$ be a smooth
family of Calabi-Yau schemes with a maximal degeneracy point at
$a_{\mathbb{C}}\in \overline C_{\mathbb{C}}$ and let
  $\overline \pi_{\mathbb{C}}: \overline X_{\mathbb{C}} \to \overline C_{\mathbb{C}}$
be a semi-stable morphism which extends $\pi_{\mathbb{C}}$. Set
$$\overline {\cal H}= H^n_{log}(\overline X_{\mathbb{C}}/\overline C_{\mathbb{C}})\, , \,
{\cal H}_{\mathbb{Z}}= Im( R^n \pi_{\mathbb{C} *} \mathbb{Z} \to R^n
\pi_{\mathbb{C} *} \mathbb{Q}). $$ Denote by
$W_{\cdot}\Psi_t^{an,un}( {\cal H}_{\mathbb{Z}})\subset \Psi_t^{an,
un}( {\cal H}_{\mathbb{Z}}) $ the monodromy filtration. \footnote{By
definition, this is a unique filtration such that the quotients
$\Psi_t^{an, un}( {\cal H}_{\mathbb{Z}})/ W_{\cdot}\Psi_t^{an, un}(
{\cal H}_{\mathbb{Z}})$ are torsion free, $N_B(W_i\Psi_t^{an, un}(
{\cal H}_{\mathbb{Z}}))\subset W_{i-2}\Psi_t^{an, un}( {\cal
H}_{\mathbb{Z}})\otimes \mathbb{Q}$ and
 $N_B^i: Gr^{n+i}_W \Psi_t^{an, un}( {\cal H}_{\mathbb{Z}}) \otimes \mathbb{Q}\simeq Gr^{n-i}_W \Psi_t^{an, un}( {\cal H}_{\mathbb{Z}}) \otimes \mathbb{Q} $.}
We then consider the limit Hodge structure
$$\Psi_t^{Hodge , un}( {\cal H})= (W_{\cdot}\Psi_t^{an, un}( {\cal H}_{\mathbb{Z}})\subset
\Psi_t^{an, un}( {\cal H}_{\mathbb{Z}})\, ,\,
a_{\mathbb{C}}^*\overline{{\cal F}}^{\cdot}\subset  a_{\mathbb{C}}^*
\overline {\cal H}  ).$$ Since $\overline{{\cal F}}^{n+1}= 0$, we
have
$$W_{-1}\Psi_t^{an, un}( {\cal H}_{\mathbb{Z}}) =0 \, , \, W_{2n} \Psi_t^{an, un}( {\cal
H}_{\mathbb{Z}})= \Psi_t^{an, un}( {\cal H}_{\mathbb{Z}})$$
$$  Im
\, N^n_B = W_{0}\Psi_t^{an, un}( {\cal H}_{\mathbb{Z}})\otimes
\mathbb{Q}.$$ Furthermore, since $ rk\, \overline{{\cal F}}^n =rk \,
\overline{{\cal F}}^0/\overline{{\cal F}}^1 = 1$ and $Im\, N_B^n\ne
0$ we must have
 $$rk \, W_0 \Psi_t^{an, un}( {\cal H}_{\mathbb{Z}}) = rk \, W_1 \Psi_t^{an, un}(
{\cal H}_{\mathbb{Z}})=1.$$ It follows that the map $$N_B^{n-1}:
\Psi_t^{an, un}( {\cal H}_{\mathbb{Z}}) \to Im \, N^{n-1}_B/ Im \,
N^n_B$$ factors through the quotient  $\Psi_t^{an, un}( {\cal
H}_{\mathbb{Z}})/ W_{2n-1} \Psi_t^{an, un}( {\cal H}_{\mathbb{Z}})$
of rank 1. In particular, $dim \, Im \, N^{n-1}_B = 2$.

Note that for any monodromy invariant lattice $P_{\mathbb{Z}}
\subset \Psi_t^{an, un}( {\cal H}_{\mathbb{Z}})$ there is a unique
local system ${\cal P}_{\mathbb{Z}}\subset  {\cal H}_{\mathbb{Z}}$
over a punctured disk $D^*\subset C_{\mathbb{C}}$ such that
$\Psi_t^{an, un}({\cal P}_{\mathbb{Z}})=P_{\mathbb{Z}}$. We apply
this remark to $L_{\mathbb{Z}}= Im \, N^{n-1}_B \cap \Psi_t^{an,
un}( {\cal H}_{\mathbb{Z}})$ and to $W_{\cdot} \Psi_t^{an,un}( {\cal
H}_{\mathbb{Z}})$. Call the corresponding local systems by ${\cal
L}_{\mathbb{Z}}$ and ${\cal W}_{\cdot} {\cal H}_{\mathbb{Z}}$. We
claim that
$${\cal L}^{Hodge}= ( {\cal W}_{\cdot} {\cal L}_{\mathbb{Z}} \, , \, {\cal F}^{\cdot}{\cal L}),$$
where ${\cal W}_{-1} {\cal L}_{\mathbb{Z}}=0$, ${\cal W}_0 {\cal
L}_{\mathbb{Z}}
 = {\cal W}_1 {\cal L}_{\mathbb{Z}}=
 {\cal W}_{0} {\cal H}_{\mathbb{Z}} $, ${\cal W}_2 {\cal
L}_{\mathbb{Z}} = {\cal L}_{ \mathbb{Z}}$ and ${\cal F}^{2}{\cal
L}=0$,  ${\cal F}^{1}{\cal L} =
 {\cal F}^{1}\cap {\cal L}$, ${\cal F}^{0}{\cal L} = {\cal
L}= {\cal L}_{ \mathbb{Z}} \otimes {\cal O}_{D^*}$, is  an
admissible variation of mixed Hodge structure
 over a sufficiently small
punctured disk $D^*$.\footnote{Recall that a variation of
 mixed Hodge structure $( {\cal W}_{\cdot} {\cal L}_{\mathbb{Z}} \subset {\cal L}_{ \mathbb{Z}}\, , \,
 {\cal F}^{\cdot}{\cal L} \subset  {\cal L})$
 over $D^*$ is called admissible if the Hodge filtration ${\cal
 F}^{\cdot}$ extends to the Deligne extension $\overline{\cal L}$ of
 ${\cal L}$.} Indeed,  over a small disk
 we have ${\cal W}_{0} \overline {\cal H} \oplus \overline{\cal F}^1 =
 \overline {\cal H} $. The claim follows.

Thus we get a class
$$[{\cal L}^{Hodge}]\in Ext^1_{VMHS}({\cal L}^{Hodge}/{\cal W}_{0}{\cal
L}^{Hodge}, {\cal W}_{0}{\cal L}^{Hodge} )$$ $$\simeq  Ext^1_{VMHS}(
\mathbb{Z}(-1), \mathbb{Z}(0)) \otimes Hom_{\mathbb{Z}} ( L_{\mathbb
{Z}}/ W_{0} \Psi_t^{an,un}( {\cal H}_{\mathbb{Z}}),     W_{0}
\Psi_t^{an,un}( {\cal H}_{\mathbb{Z}})).$$ Define
\begin{equation}\label{ccc}
q_{\mathbb C} \in Ext^1_{VMHS}(  \mathbb{Z}(-1),  \mathbb{Z}(0))
\otimes \mathbb{Q}
\end{equation}
to be the composition of $[{\cal L}^{Hodge}]$ with $$N_B^{-1}:W_{0}
\Psi_t^{an,un}( {\cal H}_{\mathbb{Z}}) \otimes \mathbb{Q} \to
  L_{\mathbb{Z}}/  W_{0} \Psi_t^{an,un}( {\cal H}_{\mathbb{Z}})   \otimes \mathbb{Q}.$$
 If the monodromy is small  i.e. that the map
$N_B: L_{\mathbb{Z}}/  W_{0} \Psi_t^{an,un}( {\cal H}_{\mathbb{Z}})
\to W_{0} \Psi_t^{an,un}( {\cal H}_{\mathbb{Z}}) $ is an
isomorphism, the class $q_{\mathbb C}$ lifts canonically to
 $$\tilde q_{\mathbb C}\in  Ext^1_{VMHS}( \mathbb{Z}(-1),
\mathbb{Z}(0)).$$  Recall that the group $Ext^1_{VMHS}(
\mathbb{Z}(-1),  \mathbb{Z}(0))$ of admissible extensions is
canonically identified with the group  of invertible functions on
$D^*$ with a regular singularity at the origin. The following lemma
immediately follows from the construction.
\begin{lemma}
The class  $\tilde q_{\mathbb{C}}$  is equal to the canonical
coordinate $q$.
\end{lemma}
We shall compute the logarithmic derivative
$$q_{\mathbb{C}}\in (\mathbb{C}((t)))^* \otimes \mathbb{Q}
\stackrel{d\, log}{\longrightarrow} \mathbb{C}[[t]] \frac {dt}{t}.$$
Let $e^0$ be a nonzero parallel section of ${\cal W}_0\overline
{\cal L}$. Then there exists a unique section  $e^1$  of $
  {\cal F}^{1}\overline{\cal L}$ such that the projection of $e^1$ to
  $\overline{\cal L}/{\cal W}_0 \overline{\cal L}$ is parallel and
$$- \frac{1}{2\pi i} N_B(e^1( a_{\mathbb{C}}))= N_{DR}(e^1(
a_{\mathbb{C}}))= e^0( a_{\mathbb{C}}).$$
 We have then
\begin{equation}\label{deq}
\nabla e^1 =  e^0 \otimes d\, log q_{\mathbb{C}}   .
\end{equation}
Assume that we are in the situation of 1.3. It follows then that
$$log q_{\mathbb{C}}\in (1+ t \mathbb{Q}[[t]] ) \frac{dt}{t}. $$
Indeed, we can normalize  $e^0$ such that $e^0( a_{\mathbb{Q}})\in
a_{\mathbb{Q}}^*H^n_{log}(\overline X_{\mathbb{Q}}/\overline
C_{\mathbb{Q}})$. Then $e^0 , e^1 \in H^n_{log}(\overline
X_{\mathbb{Q}}/\overline C_{\mathbb{Q}})$, and we are done by
(\ref{deq}).

\section{p-adic Hodge Theory}
{\bf 3.1. Fontaine-Laffaille modules.}  Fontaine-Laffaile modules
over a scheme  is a p-adic analog of variations of Hodge structure.
Below we recall this notion in the special case of torsion free
modules over a punctured disk. This is sufficient for our
applications. \footnote{Except for the last section where we need
the category of all Fontaine-Laffaille modules over a point.} The
general definition can be found in [Fa].

 Let $(D, a)$ be a formal disk over $\mathbb{Z}_p$ with a point $a: spec\,
\mathbb{Z}_p \hookrightarrow D$. We view $(D, a)$ as a logarithmic
scheme (see [Il]). A logarithmic morphism $G: (D,a) \to (D', a')$ is
morphism such that the scheme theoretical preimage of the section
$a'$ is supported on $a$ i.e. $G^*(t)= t^{\prime n} f(t')$, where
$t$ and $t'$ are coordinates on $D$ and $D'$ respectively, such that
$t(a)= t'(a')=0$, and $f$ is an invertible function on $D'$. Denote
by  $\Omega ^1 (log)$ the space of 1-forms on $D$ with logarithmic
singularities at $a$.

Let $\overline{\cal E}$ be a vector bundle over $D$ with a
logarithmic connection $\nabla : \overline{\cal E} \to
\overline{\cal E} \otimes \Omega ^1 (log)$. For any logarithmic
morphism $G: (D', a') \to (D, a)$,  $G^*\overline{\cal E}$ is
endowed with the induced logarithmic connection. Moreover, if two
logarithmic morphisms $G$ and $G'$ are equal modulo $p$ we have a
canonical parallel isomorphism
\begin{equation}\label{crystal}
\theta: G^*\overline{\cal E} \simeq G^{\prime  *}\overline{\cal E}.
\end{equation}
 In coordinates  $\theta: \overline{\cal E}\otimes _G  \mathbb{Z}_p[[t']] \to   \overline{\cal E}\otimes _{G'}  \mathbb{Z}_p[[t']] $
 is given by Taylor's formula
$$\theta (e\otimes 1) = \sum_{i=0}^{\infty} (\nabla_{\delta})^i  e
\otimes  \frac{(log(G^*(t)/G^{\prime *}(t)))^i}{i!},$$ where $\delta
= td/dt$ is the vector field.  One readily checks that
$\frac{(log(G^*(t)/G^{\prime *}(t)))^i}{i!} \in \mathbb{Z}_p[[t]] $
and that the series converges.


Let  $\tilde F: (D, a)\to (D, a)$  be a logarithmic lifting of the
Frobenius morphism (i.e. $\tilde F^*(t)= t^p (1 +p h(t) )$, where
$h(t) \in \mathbb{Z}_p[[t]]$ ). A (torsion free) Fontaine-Laffaille
module over the logarithmic disk $(D, a)$ amounts to the following data:\\
i) a vector bundle $\overline{\cal E}$ over $D$ with a filtration by
sub-bundles
$$0= \overline{\cal F}^{p-1}\subset \overline {\cal F}^{p-2} \subset \cdots \subset \overline {\cal F}^0 = \overline{\cal E},$$
ii) a logarithmic connection $\nabla : \overline{\cal E} \to
\overline{\cal E} \otimes \Omega ^1 (log)$ satisfying the Griffiths
transversality condition:
$\nabla(\overline {\cal F}^i)\subset \overline {\cal F}^{i-1} \otimes \Omega ^1 (log)$,\\
iii) a parallel morphism (``Frobenius'')
$$\phi: \tilde F^*\overline{\cal E} \to \overline{\cal E}$$
with the following properties $$\phi(\tilde F^* \overline {\cal
F}^i)\subset p^i \overline{\cal E}$$ and
 $$\sum _i p^{-i}\phi (\tilde F^*\overline {\cal F}^i)= \overline{\cal E}.$$

{\bf Remark.} The definition we gave above depends on the choice of
a lifting $\tilde F$.  Still, the categories corresponding to
different liftings are canonically equivalent. To see this  let
$\tilde F'$ be another logarithmic lifting. By (\ref{crystal}) there
is a canonical parallel isomorphism $$ \theta: \tilde F^*
\overline{\cal E} \simeq \tilde F^{\prime *}\overline{\cal E}. $$
 The functor, that provides the equivalence, takes $(\overline{\cal E}, \overline {\cal F}^i ,  \nabla) $
  to the same objects and sends $\phi$ to $\phi \theta ^{-1}$.
The  Griffiths transversality condition implies  that
$$ \theta ^{-1} \tilde F^{\prime *} \overline {\cal F}^k
\subset \sum_{i\geq 0} \frac{p^i}{i!} \tilde F^* \overline {\cal
F}^{k-i} \subset \tilde F^* \overline{\cal E}.$$ Thus,  thanks to
the assumption on the range of the Hodge filtration ( $\overline
{\cal F}^{p-1}=0$ and $\overline {\cal F}^0=\overline{\cal E}$),
 $(\overline{\cal E}, \overline {\cal F}^i ,  \nabla,   \phi \theta ^{-1}) $ satisfies the requirements in iii).
Denote the category of Fontaine-Laffaille modules over $(D, a)$ by
$MF_{[0,p-2]}(D^*)$.
 A similar construction is used to define  the pullback functor
\begin{equation}\label{pullback}
G^*:MF_{[0,p-2]}(D^*)\to MF_{[0,p-2]}(D^{'*}),
\end{equation}
 for a logarithmic
morphism $G: (D', a') \to (D, a)$.

{\bf 3.2. Dwork's Lemma.} Denote by $\mathbb{Z}_p(-k)$, ($k\geq 0$),
the constant variation: $\overline{\cal E}= {\cal O}$,
$\overline{\cal F}^k=\overline{\cal E}$, $\overline{\cal
F}^{k+1}=0$, $\phi = p^k Id$. Let $ {\cal O}(D^*)$  be the space of
functions on the punctured disk (i.e. $ {\cal O}(D^*)=
\mathbb{Z}_p((t))$).
\begin{lemma}
The group $Ext^1_{MF_{[0,p-2]}(D^*)}(   \mathbb{Z}_p(-1),
\mathbb{Z}_p(0))$ is canonically isomorphic (i.e. the isomorphism
does not depend on the lifting $\tilde F$) to p-adic completion of
the group $   {\cal O}^*(D^*)$:
\begin{equation}\label{ext}
\hat {\cal O}^*(D^*):=\underset  \longleftarrow  \lim  {\cal
O}^*(D^*)/({\cal O}^*(D^*)) ^{p^i} \stackrel{\sim}{\to} Ext^1(
\mathbb{Z}_p(-1),  \mathbb{Z}_p(0)).
\end{equation}
\end{lemma}
 \begin{proof}
 Let $t$ be a coordinate on $D$, and let  $\tilde F:D\to D $ send $t$ to $t^p$.
 Consider an extension $(\overline{\cal E}, \overline {\cal F}^i, \phi)$:
$$0\to \mathbb{Z}_p(0) \to  \overline{\cal E} \to    \mathbb{Z}_p(-1) \to 0$$
Note that $\overline {\cal F}^0=\overline{\cal E}, \overline {\cal
F}^2=0$, and that last map in the exact sequence defines an
isomorphism $\overline {\cal F}^1\simeq {\cal O}_{D}$. Let $e_1 \in
\overline {\cal F}^1$ be the preimage of $1$ under the above
isomorphism, and let $e_0\in \overline{\cal E}$ be the image of $1$
under the first map in the exact sequence. Then $e_0, e_1$ form a
basis for $\overline{\cal E}$. We have:
$$\nabla e_0=0, \nabla e_1 = e_0 \otimes \omega $$
$$\phi (\tilde F^* (e_0)) = e_0, \phi( \tilde F^*( e_1))= p e_1 + p h e_0  ,$$
for some $\omega \in  \Omega ^1 (log )$ and $h\in {\cal O}(D)$.
Since $\phi$ is parallel, $\phi \nabla = \nabla \phi$. This amounts
to the following equation
$$1/p\tilde F^*\omega - \omega = dh .$$
Thus the set of extensions is in a bijection with the set of pairs
$(\omega, h)$ satisfying the above equation. One can easily see that
the above bijection is compatible with the group structure
\footnote{The group structure on the set of pairs $(\omega, h)$ is
defined by the formula $(\omega, h)+ (\omega ', h') =   (\omega +
\omega ', h+ h')  $.}  Define a homomorphism
\begin{equation}\label{reg}
 {\cal O}^*(D^*)    \to Ext^1_{MF_{[0,p-2]}(D^*)}(   \mathbb{Z}_p(-1),  \mathbb{Z}_p(0)).
\end{equation}
sending an invertible function $q$ to the pair $( d\,log q,
\frac{1}{p} log \frac{\tilde F^*q}{q^p})$. One readily sees that
(\ref{reg}) extends to the p-adic completion of ${\cal O}^*(D^*)$.
This is the map in   (\ref{ext}).  The injectivity of (\ref{ext}) is
clear from the definition and the surjectivity
 is the content of  the Dwork's lemma \footnote{
  The  Dwork's lemma is the following statement:

   Let $\omega \in  \Omega ^1 _{log }= \mathbb{Z}_p[[t]] \frac{dt}{t}$ with $Res_{0} \omega \in \mathbb{Z}$.
   The following two conditions are equivalent:

i)$1/p\tilde F^*\omega - \omega = dh $, for some $h\in
\mathbb{Z}_p[[t]]$

ii) $\omega = d\,log q$, for some $q\in {\cal O}^*(D^*)  $ .}. Let
us check that  (\ref{ext}) is independent of the choice of the
coordinate. Indeed, let $t'$ be another coordinate and let $\tilde
F'$ be the corresponding lifting of the Frobenius. The isomorphism
(\ref{crystal}):
$$\theta: \tilde F ^*\overline{\cal E} \simeq \tilde F ^{\prime  *}\overline{\cal E} $$
takes $\tilde F^*( e_0)$ to $\tilde F^{\prime *}( e_0)$ and $\tilde
F^* (e_1)$ to
$$\sum_{i=0}^{\infty}\frac{(log(\tilde F^*(t)/\tilde F^{\prime
*}(t)))^i}{i!}\,
  \tilde F^{\prime *}((\nabla_{\delta})^i  e_1)
 =\tilde F^{\prime *} e_1 + log \frac{\tilde F^* (q)}{\tilde F^{\prime
 *}(q)}\,  e_0 .  \footnote{The last equality follows from the multiplicative version of  Taylor's formula
 $f(e^b a) = (exp(b \delta )(f))(a)= f(a)+ \delta f (a) b + \frac{\delta ^2 f(a)}{2!} b^2 +\cdots
 $.}
  $$
The claim follows.
\end{proof}

{\bf 3.3. Limit Fontaine-Laffaille module.} Let $t$ be a coordinate,
$\tilde F$ the corresponding lifting of the Frobenius, and let
$(\overline{\cal E}, \overline{\cal F}^{\cdot}, \nabla , \phi) \in
MF_{[0, p-2]}(D^*)$ be a Fontaine-Laffaille module. We define
$$\Psi^{FL}_t((\overline{\cal E}, \overline{\cal F}^{\cdot}, \nabla , \phi))=
(E= a^*\overline{\cal E}, F^{\cdot}= a^*\overline{\cal F}^{\cdot},
\phi_a).$$ $\Psi^{FL}_t((\overline{\cal E}, \overline{\cal
F}^{\cdot}, \nabla , \phi))$ is a Fontaine-Laffaille module over the
point. The residue of $\nabla$ is a morphism of Fontaine-Laffaille
modules:
$$N_{DR}= Res \, \nabla : \Psi^{FL}_t((\overline{\cal E}, \overline{\cal F}^{\cdot}, \nabla ,
\phi))\to \Psi^{FL}_t((\overline{\cal E}, \overline{\cal F}^{\cdot},
\nabla , \phi))(-1).$$ In particular,
\begin{equation}\label{rel}
  N_{DR} \phi_a =p \phi_a N_{DR}.
 \end{equation}

 {\bf Remark.} The functor $
\Psi^{FL}_t$
 depends on the choice of a coordinate $t$. If $t'= bt +\cdots$ , $b\in \mathbb{Z}_p^*$ is  another coordinate
we have
$$ \Psi^{FL}_{t'}((\overline{\cal E}, \overline{\cal F}^{\cdot}, \nabla ,
\phi \theta^{-1}))\simeq (E, F^{\cdot}, \phi_a exp(N_{DR} log\,
b^{p-1} )).$$ In particular, $ \Psi^{FL}_t$ does not get changed if
we replace $t$ by $t'$ with the same derivative.

{\bf 3.4. p-adic canonical coordinate.} Let $\overline{\pi}:
\overline X\to D$ be a proper semi-stable morphism. For any $k<p-1$,
Faltings constructed in [Fa] a Fontaine-Laffaille structure on the
logarithmic De Rham cohomology $H_{log}^k(\overline X/D)$. In the
rest of this section we assume that  $n:= dim_ D \overline X <p-1 $,
and let $\overline{\cal E}: = H_{log}^{n} (\overline X/D)/
\textstyle{p-torsion} \in MF_{[0, p-2]}(D^*)$.
 The cup product  $H_{log}^{n}(\overline X/D) \otimes H_{log}^{n}(\overline X/D)  \to  H_{log}^{2 n}(\overline X/D)
 \simeq  \mathbb{Z}_p(-n) $
  induces a perfect paring
\begin{equation}\label{paring}
<\cdot, \cdot>_{DR}: \overline{\cal E}\otimes \overline{\cal E} \to
\mathbb{Z}_p(-n).
\end{equation}
 In particular,
\begin{equation}\label{reltwo}
<\phi(\tilde F^* v) , \phi (\tilde F^*u)>_{DR}= p^n< v , u>_{DR}.
\end{equation}
 Assume that $X$ is  a Calabi-Yau  scheme over $D^*$ and that $a$
is the maximal degeneracy point. That means, by definition, that
$dim \, F^n \otimes \mathbb{Q} =1$ and the operator $N_{DR}^n: E\to
E$ is not equal to $0$. Assume, in addition, that $\overline{\pi}:
\overline X \to D$ extends to a semi-stable scheme over a curve. We
have then
\begin{equation}\label{pmor}
   rk\, Im\, N_{DR}^n =1 \, , \,  rk \, Im \, N_{DR}^{n-1} =2.
\end{equation}
   This follows from Lemma \ref{M} and "the Lefschetz principle".

\begin{lemma}\label{Fr} The Frobenius operator $\phi_a$ restricted to $Im\, N_{DR}^n$ is equal to $\pm Id$.
\end{lemma}
\begin{proof} The lemma follows immediately from (\ref{rel}) and
(\ref{reltwo}).
\end{proof}

The above lemma implies the existence of a parallel section of
$\overline{\cal E}$. Namely we have the following result.
 \begin{lemma}\label{fs} Let $\overline{\cal E}$ be a vector bundle over $D$ with a logarithmic connection and
 $\phi: \tilde F ^* \overline{\cal E} \to \overline{\cal E}$ be a parallel morphism.
For any  element $w\in E$ such that $\phi_a (w)=\pm w$ there exists
a unique parallel section $s $ of $\overline{\cal E}$ with $s(a)=w$.
The section $s$ satisfies the property $\phi(\tilde F^* s) =\pm s$.
\end{lemma}
 \begin{proof}
The uniqueness part is clear. To prove  the existence we start with
any section $s'$ of $\overline{\cal E}$ with $s'(a)=w$ and consider
the sections $s'_k= (\phi \tilde F^*)^{2k}(s')$. It is easy to see
that $\nabla s'_k \in p^{2k} \overline{\cal E} \otimes \Omega ^1
(log)$ and that $s'_{k}(a)= w$. This implies that the limit
$$ s= \underset  \longrightarrow \lim s'_{k} $$
exists and satisfies all the required properties.
\end{proof}

The nilpotent operator $N_{DR}: E\to E $ gives rise to a canonical
filtration $W_0= W_1 \subset W_2 \subset \cdots \subset W_{2n}=E$ by
Fontaine-Laffaille submodules. It is a unique filtration with
torsion free quotients $W_{i+1}/W_i$ such that $W_{\cdot}\otimes
\mathbb{Q}_p$ is the monodromy filtration on $E\otimes
\mathbb{Q}_p$. The Frobenius $\phi$ preserves the filtration
$W_{\cdot}$ and $N_{DR}(W_i)\subset W_{i-2}$. Let $L^{FL}: = Im\,
N_{DR}^{n-1} \otimes \mathbb{Q}_p\cap E$. This is a
Fontaine-Laffaille submodule of $E$. It follows from (\ref{Fr}) that
the eigenvalue of $\phi$ on $W_0$ (resp. $L^{FL}/W_0$ ) is equal to
$\pm 1$ (resp. $\pm p$). Lemma (\ref{fs}) implies that the inclusion
$W_0\hookrightarrow E$ extends uniquely to a parallel morphism
${\cal W}_0:= W_0 \otimes _{\mathbb{Z}_p} {\cal O}_D \hookrightarrow
\overline{\cal E}$. Note that the projection ${\cal W}_0
\hookrightarrow \overline{\cal E} \to \overline{\cal
F}^0/\overline{\cal F}^1$ is an isomorphism. Thus the Frobenius
$\phi :\tilde F^* (\overline{\cal E}/{\cal W}_0) \to \overline{\cal
E}/{\cal W}_0 $ is divisible by $p$. Applying (\ref{fs}) again to
$\frac{\phi}{p}: \tilde F^*(\overline{\cal E}/{\cal W}_0) \to
\overline{\cal E}/{\cal W}_0$ we conclude that the inclusion
$L^{FL}/ W_0\hookrightarrow E/W_0$ extends uniquely to a parallel
morphism $L^{FL} / W_0 \otimes _{\mathbb{Z}_p }{\cal O}_D
\hookrightarrow \overline{\cal E}/ {\cal W}_0$. Finally, let ${\cal
L}^{FL}\subset \overline{\cal E}$ be the preimage of $L^{FL}/ W_0
\otimes _{\mathbb{Z}_p} {\cal O}_D$ in $\overline{\cal E}$. By
construction, ${\cal L}^{FL}$ is a unique Fontaine-Laffaille
submodule of $\overline{\cal E}$ with $\Psi^{FL}_t({\cal L}^{FL})=
L^{FL}$.
Thus we get a canonical class
$$[{\cal L}^{FL}]\in
Ext^1_{MF_{[0, p-2]}(D^*)}( {\cal L}^{FL}/{\cal W}_0 , {\cal
W}_0)\simeq $$
$$Ext^1_{MF_{[0, p-2]}(D^*)}( \mathbb{Z}_p(-1) ,
\mathbb{Z}_p(0))\otimes  Hom_{\mathbb{Z}_p}(L^{FL}/W_0, W_0).
$$
Composing this with $N_{DR}^{-1}\in Hom( W_0 \otimes \mathbb{Q}_p,
L^{FL}/ W_0 \otimes \mathbb{Q}_p )$ we get the "p-adic canonical
coordinate":
\begin{equation}\label{pcc}
q_{\mathbb{Z}_p}\in  Ext^1_{MF_{[0, p-2]}(D^*)}( \mathbb{Z}_p(-1) ,
\mathbb{Z}_p(0))  \otimes _{\mathbb{Z}_p} {\mathbb{Q}_p}\simeq \hat
{\cal O}^*(D^*) \otimes _{\mathbb{Z}_p} {\mathbb{Q}_p}.
\end{equation}
Observe that the order
$$ ord:\,  \hat {\cal O}^*(D^*) \otimes _{\mathbb{Z}_p} {\mathbb{Q}_p} \to \mathbb{Q}_p $$
of $q_{\mathbb{Z}_p}$ is equal to $1$. In particular,
$q_{\mathbb{Z}_p}\in {\cal O}^*(D) \otimes _{\mathbb{Z}}
{\mathbb{Q}}$.

Let $e^0$ be a nonzero parallel section of ${\cal W}_0\otimes
\mathbb{Q}_p$ and let $e^1$ be a section of $(\overline{\cal F}^1
\cap {\cal L}^{FL})\otimes \mathbb{Q}_p $ whose projection to
$({\cal L}^{FL}/{\cal W}_0)\otimes \mathbb{Q}_p$ is parallel and
such that $N_{DR}(e^1(a))= e^0(a)$. We then have
\begin{equation}\label{pdeq}
\nabla e^1 =  e^0 \otimes d\, log q_{\mathbb{Z}_p}   .
\end{equation}

We shall {\it the p-adic monodromy is small}, if the operator
$N_{DR}:  L^{FL}/W_0 \to  W_0 $ is an isomorphism. If this is the
case, one has
\begin{equation}\label{pint}
q_{\mathbb{Z}_p}\in {\cal O}^*(D^*)/\mu_{p-1} \subset \hat {\cal
O}^*(D^*) \otimes \mathbb{Q}_p .
\end{equation}

{\bf 3.5. p-adic Yukawa map.} In this subsection we assume that the
p-adic monodromy is small. Denote by $q\in {\cal O}_{D}$ the p-adic
canonical coordinate (defined up to a (p-1)th root of unity). Let
\begin{equation}\label{ks}
S^{n}T_{D, log} \to Hom _{{\cal O}_D} (\overline{\cal F}^n,
\overline{\cal F}^0/\overline{\cal F}^1)\simeq (\overline{\cal F}^n
\otimes \overline{\cal F}^n)^*
\end{equation}
be the Kodaira-Spenser morphism. Here $T_{D, log}$ denotes the sheaf
dual to $\Omega^1(log)$ i.e. the sheaf of vector fields on $D$
vanishing at $a$. Choose a generator $e^0 $ of ${\cal W}_0^{\nabla}$
and let $e_0\in \overline{\cal F}^n$ be a section with $(e^0, e_0)=
1$ \footnote{The paring ${\cal W}_0\otimes \overline{\cal F}^n \to
{\cal O}_D$ is perfect. For the projection ${\cal W}_0
\hookrightarrow \overline{\cal E} \to \overline{\cal
F}^0/\overline{\cal F}^1$ is an isomorphism.}. Applying (\ref{ks})
to $(q \frac{d}{dq})^{\otimes ^n} $ and pairing the result with
$e_0\otimes e_0 $ we obtain the p-adic Yukawa function $Y \in {\cal
O}_D$. Observe that $Y(q)$ is well defined up to multiplication by a
constant in $\mathbb{Z}_p^*$.

\begin{proposition}\label{kont} Assume that $n=3$ and that $rk\, E = 4$. Then
$$Y(q)= n_0 + \sum^{\infty}_{d=1}n_d d^3 \frac{q^d}{1-q^d},$$
where $n_d \in \mathbb{Z}_p$.
\end{proposition}
\begin{proof}
We shall use the following elementary result:
\begin{lemma}\label{KSV} ([KSV]. Lemma 2.)
Assume that a formal power series  $Y(q)\in \mathbb{Z}_p[[q]]$ is
written in the form
$$Y(q)=    \sum^{\infty}_{d=1}n_d d^3 \frac{q^d}{1-q^d}.$$
Then $n_d \in \mathbb{Z}_p$ if and only if $Y(q)- Y(q^p)=
\delta^3(\psi(q))$, for some $\psi(q)\in \mathbb{Z}_p[[q]]$. Here
$\delta= q\frac{d}{dq}$.
\end{lemma}
\begin{lemma}\label{ht} The monodromy filtration $W_0= W_1 \subset W_2= W_3 \subset W_4=W_5  \subset W_6= E$,
extends to a filtration ${\cal W}_i \subset \overline{\cal E} $ by
Fontaine-Laffaille submodules such that either ${\cal W}_{2i}/{\cal
W}_{2i-2} \simeq \mathbb{Z}_p(-i)$, for all $0\leq i \leq 3$, or
${\cal W}_{2i}/{\cal W}_{2i-2} \simeq \epsilon \mathbb{Z}_p(-i)$.
Here $\epsilon \mathbb{Z}_p(-i)$ denotes the constant
Fontaine-Laffaille module with
 $\overline{\cal F}^i={\cal O}_D$, $\overline{\cal F}^{i+1}=0$, $\phi = -p^i Id$.
\end{lemma}
\begin{proof}
Our assumptions imply that  $rk\, W_{2i}/W_{2i-2}= 1$, for $0\leq i
\leq 3$. Thus, by Lemma \ref{Fr} and (\ref{rel}), the operator
$\phi$ acts on   $W_{2i}/W_{2i-2}$ as   $\pm p^i Id$.

We prove by induction on $i$ that $W_{2i}\subset E$ extends to a
subbundle ${\cal W}_{2i}$ of $\overline{\cal E}$ preserved by the
connection and that  $\overline{\cal F}^{i+1}\oplus {\cal W}_{2i}=
\overline{\cal E}$. Indeed, for $i=-1$, there is nothing to prove.
Assume that we know the result for $i=k$. Then $(\overline{\cal
E}/{\cal W}_{2k}=(\overline{\cal F}^{k+1}+{\cal W}_{2k})/ {\cal
W}_{2k}\supset \cdots \supset (\overline{\cal F}^{3}+{\cal W}_{2k})/
{\cal W}_{2k}, \phi, \nabla) $ is a Fontaine-Laffaille module.
Applying Lemma \ref{fs} to $\overline{\cal E}/{\cal W}_{2k} \otimes
\mathbb{Z}_p(k+1)$ we see that $W_{2k+2}/W_{2k}\subset E/W_{2k}$
extends to a subbundle of ${\cal W}_{2k+2}/{\cal W}_{2k}\subset
\overline{\cal E}/{\cal W}_{2k}$. It remains to show that ${\cal
W}_{2k+2}/{\cal W}_{2k}\oplus (\overline{\cal F}^{k+2}+{\cal
W}_{2k})/ {\cal W}_{2k}= \overline{\cal E}/{\cal W}_{2k}$. We will
be done if we prove that this is true over the closed point of $D$
i.e. $(W_{2k+2}/{W}_{2k}\oplus (F^{k+2}+ W_{2k})/ {W}_{2k})\otimes
\mathbb{F}_p = (E/W_{2k})\otimes \mathbb{F}_p$. Indeed, the operator
$p^{-k-1} \phi$ induces an action on  $ (E/W_{2k})\otimes
\mathbb{F}_p$ which is $0$ on
 $ ((F^{k+2}+ W_{2k})/ {W}_{2k})\otimes \mathbb{F}_p$ and invertible on  $(W_{2k+2}/{W}_{2k})\otimes \mathbb{F}_p$. The claim follows.
\end{proof}

For the rest of the proof we assume that $\phi$ acts on $W_0$ as
$+Id$.  The other alternative is considered in a similar way.

By the definition of the canonical coordinate $q$ we can find
sections $e^0\in {\cal W}_0$, $e^1\in \overline{\cal F}^1\cap {\cal
W}_2$ such that
$$
\nabla_\delta  e^0 = 0, \, \phi e^0= e^0, \, \nabla_\delta e^1 =
e^0, \,   \phi e^1 = p e^1,
$$
and such that $e^0, e^1$ generate ${\cal W}_2$. Next, it follows
from Lemma \ref{ht} that there exist  unique $e_1 \in \overline{\cal
F}^2\cap {\cal W}_4$, $e_0 \in  \overline{\cal F}^3 $ such that
$$
(e^0, e_0)= 1, \, (e^1, e_1)= -1
$$
Observe that $e^i, e_i$ generate $\overline{\cal E}$. Thanks to the
self-duality condition (\ref{paring}) we have
$$
\nabla_{\delta} e_1 = Y(q)e^1 ,  \, \nabla_{\delta} e_0= e_1
$$
$$
\phi e_1 = p^2(e_1 + m_{23}(q) e^1  + m_{13}(q) e^0),  \,  \phi e_0=
p^3(e_0 - m_{13}(q) e^1 + m_{14}(q) e^0),
$$
where $Y(q)$ is the Yukawa function. Finally, the relation
$\nabla_{\delta} \phi = p \phi \nabla_\delta $ amounts to
$$ Y(q)- Y(q^p)= \delta (m_{23}), \, m_{23}= -\delta(m_{13}), \, \delta  (m_{14})= 2 m_{13}.$$
Thus
$$Y(q^p) - Y(q) = \frac{1}{2} \delta ^3 m_{14},$$
and we are done by Lemma \ref{KSV}.

\end{proof}

\section{Comparison}
{\bf 4.1. Plan of the proofs of Theorems \ref{th1} and \ref{th2}.}

Let $\overline \pi: \overline X_S \to \overline C_S$ be a
semi-stable morphism satisfying the conditions i) - iii) from
Section 1.3. Denote by $q_{\mathbb{C}}\in (\mathbb{C}((t)))^*
\otimes _{\mathbb{Z}} {\mathbb{Q}}$ the complex canonical coordinate
(\ref{ccc}) and by $q_{\mathbb{Z}_p} \in (\mathbb{Z}_p((t)))^*
\otimes _{\mathbb{Z}} {\mathbb{Q}}$ the p-adic one (\ref{pcc}).


\begin{proposition}\label{dif} a)  For every prime
prime $p$ such that $(p,N)=1$, we have
$$q_{\mathbb{C} }=  q_{\mathbb{Z}_p}  \in ({\mathbb{Q}}((t)))^* \otimes \mathbb{Q}.$$
b) Assume that the Betti monodromy of the family $\overline X_S \to
\overline C_S$ is small (see 1.1 ). Then, for every prime $p$ with
$(p,N)=1$, the p-adic monodromy is
also small (see 3.4 ).\\
c) Let $\omega $ be a nonvanishing section of  the line bundle
$\overline {\cal F}^n= \overline \pi_* \Omega_{\overline
X_S/\overline C_S}^n(log \, Y_S)$  over an open neighborhood of the
subscheme $a: S \hookrightarrow \overline C_S $. Then
$$(\frac{1}{(2\pi i)^n} \int
_{\delta _1} \omega)^2 \in  (\mathbb{Z}[N^{-1}][[t]])^*.$$

\end{proposition}

In the remaining part of this section we complete the proofs of
Theorems \ref{th1} and \ref{th2} assuming Proposition \ref{dif}. A
proof of the proposition (which is the hardest technical part of the
argument) is given in Sections 4.2-4.5.

{\bf Proof of Theorem \ref{th1}.} Since $q'(0)\in \mathbb{Q}^*$ and
$ d\, log \, q (t) \in  \mathbb{Q}[[t]] \frac{dt}{t}$ the
coefficients of $q(t)$ are rational numbers. On the other hand,
parts a) and b) of Proposition \ref{dif} together with formula
(\ref{pint}) show that, for every prime $p$ such that $(p,N)=1$,
$$q(t)\in (\mathbb{Z}_p((t)))^* \cap (\mathbb{Q}((t)))^* \subset
 (\mathbb{Q}_p((t)))^*.$$
This completes the proof.

{\bf Proof of Theorem \ref{th2}.} Let $\omega \times \omega$ be a
local  section of  $\overline \pi_* \Omega_{\overline
X_{\mathbb{C}}/\overline C_{\mathbb{C}}}^n(log \, Y_{\mathbb{C}})
\otimes \overline \pi_* \Omega_{\overline X_{\mathbb{C}}/\overline
C_{\mathbb{C}}}^n(log \, Y_{\mathbb{C}})$ defined by the equation
$$\frac{1}{(2\pi i)^n} \int
_{\delta _1} \omega =1.$$ Part c) of Proposition \ref{dif} shows
that $\omega \times \omega$ yields a nonvanishing  section of $
\overline \pi_* \Omega_{\overline X_S/\overline C_S}^n(log \,
Y_S)\otimes \overline \pi_* \Omega_{\overline X_S/\overline
C_S}^n(log \, Y_S)$ over the formal neighborhood $D_S$.
 This together with Theorem \ref{th1} imply that the coefficients of the Yukawa function
  $Y(q)$ are rational numbers and so are the instanton numbers $n_d$. It also follows that $Y(q)$ coincides
(up to a constant factor in $\mathbb{Z}_p^*$) with the p-adic Yukawa
function from {\bf 3.5}. Thus by Proposition \ref{kont} the numbers
$n_d$ are p-adic integers. This completes the proof.

{\bf 4.2. Recollections on p-adic Comparison Theorem.} Recall from
[FL] that there is an exact tensor fully faithful functor
$$U: MF_{[0,p-2]} \to Rep(\Gamma)$$
from the category $MF_{[0,p-2]}$ of Fontaine-Laffaille modules over
$spec \, \mathbb{Z}_p$  to the category $Rep(\Gamma)$ of finitely
generated $\mathbb{Z}_p$-modules equipped with an action of the
Galois group $\Gamma =Gal(\overline{\mathbb{Q}}_p / \mathbb{Q}_p)$.
We will use the following properties of $U$:

1) $U$ takes a finite (as a plain abelian group) Fontaine-Laffaille
module to a $\Gamma $-module of the same finite order.

2) $U(\mathbb{Z}_p(i))= \mathbb{Z}_p(i)$ and the induced morphism
$$ \underset  \longleftarrow  \lim \mathbb{Z}_p ^* /( \mathbb{Z}_p ^* ) ^{p^i} \simeq
 Ext^1_{MF_{[0,p-2]}}(\mathbb{Z}_p(-1), \mathbb{Z}_p(0)) \stackrel{U}{\hookrightarrow}
   Ext^1_{Rep(\Gamma)}(\mathbb{Z}_p(-1), \mathbb{Z}_p(0))$$
  $$ \stackrel{Kummer}{\simeq}
 \underset  \longleftarrow  \lim \mathbb{Q}_p ^* /( \mathbb{Q}_p ^* ) ^{p^i}$$
is identity.

3)Let $\overline {\pi}: \overline X_{\mathbb{Z}_p} \to \overline
C_{\mathbb{Z}_p} $ be a proper semi-stable (relative to
$\mathbb{Z}_p$ ) scheme.
 Assume that
$dim _{ \overline C_{\mathbb{Z}_p}} \overline X_{\mathbb{Z}_p}
\leq p-2$.  Then there is a canonical isomorphism:
\begin{equation}
 U(\Psi^{FL}_t(H^k_{log}( \overline X_{\mathbb{Z}_p}/\overline C_{\mathbb{Z}_p})) ) \stackrel{\simeq}{\longrightarrow}
   \Psi_{t}^{et} (R^k \pi_{\mathbb{Q}_p *}^{et} \mathbb{Z}_p)
\end{equation}
 Here $\pi_{\mathbb{Q}_p}$ denotes the projection $X_{\mathbb{Q}_p} \to C_{\mathbb{Q}_p}$ and
$\Psi^{et}_{t}: Sh^{et}(C_{\mathbb{Q}_p}) \to
Sh^{et}(a_{\mathbb{Q}_p})=  Rep(\Gamma) $ is the etale vanishing
cycles functor.
Moreover, we have the following commutative diagram

   $$
\def\normalbaselines{\baselineskip20pt
\lineskip3pt  \lineskiplimit3pt}
\def\mapright#1{\smash{
\mathop{\to}\limits^{#1}}}
\def\mapdown#1{\Big\downarrow\rlap
{$\vcenter{\hbox{$\scriptstyle#1$}}$}}
\begin{matrix}
 U(\Psi^{FL}_t(H^k_{log}( \overline X_{\mathbb{Z}_p}/\overline C_{\mathbb{Z}_p})) ) &    \stackrel{\simeq}{\longrightarrow}
 &  \Psi_{t}^{et} (R^k \pi_{\mathbb{Q}_p *}^{et} \mathbb{Z}_p)  \cr
\mapdown{N_{DR}}  &  &\mapdown{N_{et}}   \cr
 U(\Psi^{FL}_{t}(H^k_{log}( \overline X_{\mathbb{Z}_p}/\overline C_{\mathbb{Z}_p}))_a )\otimes \mathbb{Z}_p(-1) &    \stackrel{\simeq}{\longrightarrow}
 &  \Psi_{t}^{et} (R^k \pi_{\mathbb{Q}_p *}^{et} \mathbb{Z}_p) \otimes \mathbb{Z}_p(-1)
 \end{matrix}
 $$
This follows from the main Comparison Theorem in [Fa].

{\bf 4.3. 1-motives, the motivic Albanese functor $LAlb$.} The main
references here are [D3] and [BK]. Let $k$ be a field of
characteristic $0$. Fix an algebraic closure $\overline k \supset
k$. A 1-motive over $k$  is a triple
$$M = (\Lambda, G,  \Lambda \stackrel{u}{\longrightarrow} G(\overline{k})),$$
 where $\Lambda$ is a free abelian group of finite rank equipped with an action of the Galois group $Gal(\overline k/k)$
that factors through a finite quotient, $G$ is an semi-abelian
variety over $k$ i.e. an extension
\begin{equation}
0\to T \to G \to A\to 0
\end{equation}
of an abelian variety by a torus, and $u$ is a homomorphism of the
Galois modules. We shall denote by ${\cal M}_1(k)$ the additive
category of 1-motives \footnote{The Galois module $ \Lambda $ can be
viewed as a discrete group scheme over $spec \, k$. Giving a
homomorphism $\Lambda \longrightarrow G(\overline{k})$ of Galois
modules is equivalent to giving a morphism $ \underline \Lambda
\longrightarrow \underline G$ of the \'etale sheaves represented by
$\Lambda$ and $G$. This remark provides a construction of the
category ${\cal M}_1(k)$ that is independent of the choice of an
algebraic closure $\overline k$ ([BK]).}.

Every 1-motive is equipped with a canonical (weight) filtration:
$$W_{-2}M = (0, T) \subset W_{-1}M = (0, G) \subset W_0 M= M.$$
Thus $W_0 M/W_{-1} M = (\Lambda, 0)$ and  $W_{-1} M/W_{-2} M = (0,
A)$. The category ${\cal M}_1(k; \mathbb{Q}):= {\cal M}_1(k) \otimes
\mathbb{Q} $ is abelian ([BK], Proposition 1.1.5) and any morphism
in
 ${\cal M}_1(k; \mathbb{Q})$ is strictly compatible with the weight filtration.

Set $\mathbb{Z}(0) = (\mathbb{Z}, 0)$ and $\mathbb{Z}(1) = (0,
\mathbb{G}_m)$. The same 1-motives $\mathbb{Z}(i)$ ($i=0,1$) but
viewed as objects of  ${\cal M}_1(k; \mathbb{Q})$ are denoted by
$\mathbb{Q}(i)$. We have
\begin{equation}\label{tateext}
   Ext^1_{{\cal M}_1(k)}(\mathbb{Z}(0), \mathbb{Z}(1)) \simeq k^*,  \,
 Ext^1_{{\cal M}_1(k;  \mathbb{Q})}(\mathbb{Q}(0), \mathbb{Q}(1)) \simeq k^* \otimes \mathbb{Q}.
\end{equation}
 For any prime $p$, we have
the etale realizations functors ([D3], 10.1.5):
 $$T_{\mathbb{Z}_p}^{et}: {\cal M}_1(k) \to Rep_{\mathbb{Z}_p}(Gal(\overline k/k)),$$
 $$  T_{\mathbb{Q}_p}^{et}: {\cal M}_1(k; \mathbb{Q}) \to
Rep_{\mathbb{Q}_p}(Gal(\overline k/k)).$$  We also set
$$T_{\mathbb{Q}_p}^{*et}(M)= Hom_{\mathbb{Q}_p}( T_{\mathbb{Q}_p}^{et}(M), \mathbb{Q}_p)\in
Rep_{\mathbb{Q}_p}(Gal(\overline k/k)).$$

If $k=\mathbb{C}$ the category
 ${\cal M}_1(\mathbb{C})$ is equivalent to the category of torsion free polarizable mixed Hodge structures
of type $\{(0,0),(0,-1),(-1, 0),(-1,-1)\}$ ([D3], 10.1.3):
\begin{equation}\label{DeRham}
T^{Hodge}:  {\cal M}_1(\mathbb{C})
\stackrel{\simeq}{\longrightarrow} MHS_1
\end{equation}
 For $k\subset \mathbb{C}$, $M\in  {\cal M}_1(k)$,
$T^{Hodge}(M\times _k spec \, \mathbb{C})= (W_{\cdot}\subset
V_{\mathbb{Z}}, F^{\cdot}\subset V_{\mathbb{C}})$ there is a
functorial isomorphism of $\mathbb{Z}_p$-modules
\begin{equation}\label{be}
V_{\mathbb{Z}}\otimes \mathbb{Z}_p \simeq T_{\mathbb{Z}_p}^{et}(M).
\end{equation}
Abusing notation, we shall also denote by $T^{Hodge}$ the
equivalence $${\cal M}_1(\mathbb{C};
\mathbb{Q})\stackrel{\simeq}{\longrightarrow} MHS_1^{\mathbb Q}=
MHS_1\otimes \mathbb{Q}$$ induced by (\ref{DeRham}) and the
corresponding equivalence of the derived categories
$$D^b({\cal M}_1(\mathbb{C};
\mathbb{Q}))\stackrel{\simeq}{\longrightarrow} D^b(MHS_1^{\mathbb
Q}).$$
 Let
$$D^b({\cal M}_1(k; \mathbb{Q}))  \stackrel{{\cal S}}{\hookrightarrow} DM^{eff}_{gm}(k; \mathbb{Q})$$
be embedding of the bounded derived category of 1-motives into the
triangulated category of Voevodsky motives ([O]). By [BK] ${\cal S}$
has a left adjoint functor:
$$LAlb:   DM^{eff}_{gm}(k; \mathbb{Q}) \to   D^b({\cal M}_1(k; \mathbb{Q})).$$
Denote by $MHS^{ \mathbb{Q}}$ the category
 of mixed polarizable Hodge
structures over ${\mathbb Q}$ and by $MHS_{eff}^{\mathbb{Q}}$ the
full subcategory of $MHS^{ \mathbb{Q}}$, whose objects are mixed
Hodge structures $(W_{\cdot}\subset V_{\mathbb{Q}}, F^{\cdot}\subset
V_{\mathbb{C}})$ with
 $ F^1=0$.
It is proven in [Vol]
 that embedding of the derived categories
$$\overline {\cal S}: D^b(MHS_1^{\mathbb Q}) \to  D^b(MHS_{eff}^{ \mathbb{Q}})$$
admits a $t$-exact left adjoint functor \footnote{We say that a
triangulated functor $T: D({\cal A})\to D({\cal B})$ is t-exact if
$T({\cal A})$ belongs to the essential image of ${\cal B}$ in
$D({\cal B})$.}
$$\overline {LAlb}: D^b(MHS_{eff}^{\mathbb{Q}}) \to D^b(MHS_1^{\mathbb Q})$$
and that
$$ T^{Hodge}\circ LAlb \simeq \overline{LAlb} \circ R^{Hodge}:
DM^{eff}_{gm}(\mathbb{C}; \mathbb{Q}) \to  D^b(MHS_1^{\mathbb Q}).$$
Here
$$ R^{Hodge}: DM^{eff}_{gm}(\mathbb{C} ; \mathbb{Q})   \to  D^b(MHS^{\mathbb
Q}_{eff})$$ is the homological Hodge realization functor (i.e.
$R^{Hodge}(M)= R_{Hodge}(M)^*$, where $R_{Hodge}$ is Huber's
cohomological realization ([Hu1], [Hu2]).)

{\bf 4.4. Motivic vanishing cycles.} Let $X_k
\stackrel{\pi}{\longrightarrow} C_k $ be a smooth proper scheme over
a punctured curve $ C_k \hookrightarrow \overline C_k
\stackrel{a}{\hookleftarrow} spec \, k$ over a field $k\subset
\mathbb{C}$.
Fix a local coordinate $t$ at $a$ and an integer $m\geq 0$. Denote
by $$H^m(X_{\mathbb{C}}/C_{\mathbb{C}})= ( R^m\pi_*
\mathbb{Q},F^{\cdot}\subset
H^m_{DR}(X_{\mathbb{C}}/C_{\mathbb{C}}))$$  the variation of Hodge
structure associated to the family $X_{\mathbb{C}}
\stackrel{\pi_{\mathbb{C}}}{\longrightarrow} C_{\mathbb{C}} $ and by
$ H_m(X_{\mathbb{C}}/C_{\mathbb{C}})$ the dual variation. Let
$\Psi_t^{Hodge, un}( H_m(X_{\mathbb{C}}/C_{\mathbb{C}}))$ be the
unipotent limiting mixed Hodge structure. The Hodge structure
$$\overline {LAlb} (\Psi_t^{Hodge, un}(
H_m(X_{\mathbb{C}}/C_{\mathbb{C}})))$$ can be viewed as a 1-motive
over $\mathbb{C}$. In this subsection we explain how this 1-motive
canonically descends to a 1-motive $$M_{t,m}(X_k)= (\Lambda, G,
\Lambda \stackrel{u}{\longrightarrow} G(\overline{k}))\in {\cal
M}_1(k; \mathbb{Q})  $$ over $k$. In addition, $M_{t,m}(X_k)$ comes
equipped with a ``monodromy" homomorphism $N: \Lambda \otimes
\mathbb{Q} \to T_* \otimes \mathbb{Q} $ of $Gal(\overline
k/k)$-modules. Here $T_*$ denotes the group of cocharacters of the
torus $T$: $T_*= \underline {Hom}(\mathbb{G}_m, T)(\overline k)$.
Equivalently, $N$ can be viewed as a morphism of 1-motives:
$$N: (W_0 M_{t,m}(X_k)  /W_{-1} M_{t,m}(X_k))(1) \to W_{-2} M_{t,m}(X_k).$$
The main properties of $M_{t,m}(X_k)$ are the following.

1) There is a natural isomorphism:
\begin{equation}\label{c1}
\overline {LAlb} (\Psi_t ^{Hodge, un}(
H_m(X_{\mathbb{C}}/C_{\mathbb{C}})))\simeq T^{Hodge}(
M_{t,m}(X_k))\footnote{Abusing notation, we denote by $T^{Hodge}$
the composition of functors ${\cal M}_1(k; \mathbb{Q})\to  {\cal
M}_1(\mathbb{C}; \mathbb{Q})
\stackrel{T^{Hodge}}{\longrightarrow}MHS_1^{\mathbb Q}$.}
\end{equation}
compatible with the monodromy action.

  2) There is a natural morphism
$Gal(\overline k/k)$-modules
\begin{equation}\label{c2}
\alpha:  T^{*et}_{ \mathbb{Q}_p}(M_{t,m}(X_k)) \to  \Psi_t^{et,
un}(R^m \pi_*^{et}\mathbb{Q}_p)
\end{equation}
where
$$\Psi_t^{et, un}: Sh^{et}(C_k) \to Sh^{et}(spec\, k)\to Rep_{\mathbb{Q}_p}(Gal(\overline k/k))$$
 denotes
the functor of unipotent vanishing cycles (see [B]). The morphism
$\alpha $ commutes with the monodromy action.

3) If $k'\supset k$ is any field extension, there is a natural
isomorphism
\begin{equation}\label{c3}
M_{t,m}(X_k \times _k spec\, k')\simeq  M_{t,m}(X_k)\times _k spec\,
k'
\end{equation}
compatible in the obvious way with (\ref{c2}).

4) Assume  that $k\subset \mathbb{C}$. Set
$$T^{Hodge}(  M_{t,m}(X_k))= (W_{\cdot}\subset V_{\mathbb{Q}}, F^{\cdot}\subset V_{\mathbb{C}}).$$
  The following diagram is commutative.

   $$
\def\normalbaselines{\baselineskip20pt
\lineskip3pt  \lineskiplimit3pt}
\def\mapright#1{\smash{
\mathop{\to}\limits^{#1}}}
\def\mapdown#1{\Big\downarrow\rlap
{$\vcenter{\hbox{$\scriptstyle#1$}}$}}
\begin{matrix}

   V^*_{\mathbb{Q}}\otimes \mathbb{Q}_p  &    \stackrel{(\ref{c1})}{\longrightarrow}
 &  \Psi_t^{an, un}(R^m \pi_*^{an}\mathbb{Q})\otimes \mathbb{Q}_p          \cr
\mapdown{(\ref{be})}  &  &\mapdown{}   \cr
    T^{*et}_{ \mathbb{Q}_p}(M_{t,m}(X_k))      &    \stackrel{\alpha}{\longrightarrow}
 &  \Psi_t^{et, un}(R^m \pi_*^{et}\mathbb{Q}_p)
 \end{matrix}
 $$

The above properties of $M_{t,m}(X_k)$ are sufficient for our
applications. We shall indicate a conceptual construction of
$M_{t,m}(X_k)$ based on the theory of Voevodsky's motives.
Unfortunately, the construction relies on the following general fact
that is not explained in the published literature.

Let
$$\Psi^{mot, un}_t: DM^{eff}_{gm}(\eta; \mathbb{Q}) \to DM^{eff}_{gm}(k; \mathbb{Q})$$
 the functor of (unipotent) motivic vanishing cycles from the triangulated category of motives over the
 generic point $\eta \in C_{k}$ to the category of motives over $a$, and let
$$N: \Psi^{mot, un}_t (1) \to \Psi^{mot, un}_t $$
be the monodromy operator (see [A1]).
 The fact, we will need, is that the formation $(\Psi^{mot, un}_t, N)$
  commutes with the etale and Hodge realizations [Hu1], [Hu2]:

   $$
\def\normalbaselines{\baselineskip20pt
\lineskip3pt  \lineskiplimit3pt}
\def\mapright#1{\smash{
\mathop{\to}\limits^{#1}}}
\def\mapdown#1{\Big\downarrow\rlap
{$\vcenter{\hbox{$\scriptstyle#1$}}$}}
\begin{matrix}

  DM^{eff}_{gm}(\eta; \mathbb{Q}) &    \stackrel{\Psi^{mot, un}_t }{\longrightarrow}
 &    DM^{eff}_{gm}(k; \mathbb{Q})                       \cr
\mapdown{ }  &  &\mapdown{}    \cr
 D^b( Rep_{\mathbb{Q}_p}(Gal(\overline {k(\eta)}/k)) )       &    \stackrel{ \Psi^{et, un}_t }{\longrightarrow}
 &  D^b( Rep_{\mathbb{Q}_p}(Gal(\overline {k}/k)) )
\end{matrix}
 $$
   $$
\def\normalbaselines{\baselineskip20pt
\lineskip3pt  \lineskiplimit3pt}
\def\mapright#1{\smash{
\mathop{\to}\limits^{#1}}}
\def\mapdown#1{\Big\downarrow\rlap
{$\vcenter{\hbox{$\scriptstyle#1$}}$}}
\begin{matrix}

    DM^{eff}_{gm}(\eta ; \mathbb{Q})   &    \stackrel{\Psi^{mot, un}_t }{\longrightarrow}
 &    DM^{eff}_{gm}(\mathbb{C}; \mathbb{Q})                       \cr
   \mapdown{R^{Hodge} }  &  &\mapdown{R^{Hodge}}   \cr
 D^b(VMHS_{eff}^{\mathbb{Q}}(\eta))  & \stackrel{ \Psi^{Hodge, un}_t }{\longrightarrow}
 &  D^b(MHS^{\mathbb Q}_{eff}).
\end{matrix}
 $$

Assuming this fact we construct $M_{t,m}(X_k)$ as follows. It is
proven in [BK] the fully faithful functor $D^b({\cal M}_1(k;
\mathbb{Q}))  \to DM^{eff}_{gm}(k; \mathbb{Q})$ has a left adjoint:
$$LAlb:   DM^{eff}_{gm}(k; \mathbb{Q}) \to   D^b({\cal M}_1(k; \mathbb{Q})).$$
Set
$$M_{t,m}(X_k):= H_m(LAlb \, \Psi_t^{mot, un}(
\mathbb{Q}_{tr}[X_{\eta}])).$$
 Let us just explain that
$M_{t,m}(X_k)$ has the key property 1. Indeed, by Theorem 2 from
[Vol] the Albanese functor commutes with the Hodge realization.
Thus, we have
$$T^{Hodge}(M_{t,m}(X_{\mathbb C}))\simeq H_m(\overline{LAlb} \, R^{Hodge} \Psi_t^{mot, un}(
\mathbb{Q}_{tr}[X_{\eta}]))$$
$$\simeq  \overline {LAlb}\,
\Psi_t ^{Hodge, un} H_m(X_{\mathbb{C}}/C_{\mathbb{C}}).$$

{\bf Example.} \footnote{The reader can skip this example: it will
not be used in the main text below.}
 Here we explain an elementary
construction of the motive $M_{t,1}(X_k)$ for a family $\pi: X_k \to
C_k$  of curves with a semi-stable reduction. Choose a semi-stable
model $\overline \pi: \overline X_k \to \overline C_k$, such that
all the irreducible components $Y_{\overline k, \gamma }$ of the
special fiber $Y_{\overline k}:= \overline X_k \times _{\overline
C_k} \overline k $ are smooth. Let $\Gamma$ be the free abelian
group whose generators $[\gamma]$ correspond to irreducible
components of $Y_{\overline k}$. For each singular point $y_{\mu}$
of $Y_{\overline k}$ we denote by $R_{\mu}$ the subgroup of $\wedge
^2 \Gamma$ generated by $[\gamma _1] \wedge [\gamma_2]$, where
$Y_{\overline k, \gamma_1 }$ and $Y_{\overline k, \gamma_2 }$ are
the two components meeting at $y_{\mu}$ (i.e. $R_{\mu}$ is
isomorphic to $\mathbb{Z}$ but the isomorphism depends on the
ordering of the components meeting at $y_{\mu}$).
 Define a
homomorphism
$$u: R_{\mu} \to Pic(Y_{\overline k})$$
as follows. Consider the invertible sheaf ${\cal O}(y_{\mu,
\gamma_1} - y_{\mu, \gamma_2})$ on the normalization $\tilde
Y_{\overline k} \to Y_{\overline k}$, where $y_{\mu, \gamma_1}$,
$y_{\mu, \gamma_2}$ are the preimages of $y_{\mu}$ in $\tilde
Y_{\overline k}$. We claim that ${\cal O}(y_{\mu, \gamma_1} -
y_{\mu, \gamma_2})$ canonically descends  to a line bundle
$u([\gamma _1] \wedge [\gamma_2])$ over $Y_{\overline k}$: the
descend data are trivial outside of points $y_{\mu, \gamma_1}$,
$y_{\mu, \gamma_2}$, and the identification
$${\cal O}(y_{\mu, \gamma_1} -
y_{\mu, \gamma_2})_{y_{\mu, \gamma_1}} \simeq {\cal O}(y_{\mu,
\gamma_1} - y_{\mu, \gamma_2})_{y_{\mu, \gamma_2}}$$ is given by a
canonical isomorphism
$$ T_{Y_{\overline k, \gamma_1 }, y_{\mu, \gamma_1}} \otimes T_{Y_{\overline k, \gamma_2 }, y_{\mu,
\gamma_2}} \simeq T_{\overline C_{\overline k}, a} \simeq \overline
k. \footnote{The ismorphism $T_{\overline C_{\overline k}, a} \simeq
\overline k$ is determined by the coordinate $t$.} $$
 Finally, let $\Lambda \subset \bigoplus _\mu R_\mu $ be the kernel of the degree
 map:
 $$ \bigoplus _\mu R_\mu  \stackrel{u}{\longrightarrow}  Pic(Y_{\overline
 k}) \stackrel{deg}{\longrightarrow} \bigoplus _\gamma \mathbb{Z}.$$
Then $M_{t,m}(X_k)=(\Lambda, {\bf Pic}^0(Y_{\overline
 k}), \Lambda \stackrel{u}{\to} {\bf Pic}^0(Y_{k})(\overline k)).$

 {\bf 4.5. Proof of Proposition \ref{dif}.} a) By (\ref{deq}) and
(\ref{pdeq}) we have
\begin{equation}\label{eqcor}
d\, log \,  q_{\mathbb{C}}(t) = d\, log \, q_{\mathbb{Z}_p}(t)\in
(1+ t \mathbb{Q}[[t]] ) \frac{dt}{t} .
\end{equation}
Thus, it suffices to prove that
$$q'_{\mathbb{C} }(0)=  q'_{\mathbb{Z}_p}(0)  \in \mathbb{Q}^* \otimes \mathbb{Q} .$$

Consider the 1-motive $M_{t,n}(X_{\mathbb{Q}})= (\Lambda, G, \Lambda
\stackrel{u}{\longrightarrow} G(\overline{k})) $. By (\ref{c1}) we
have
$$L^{Hodge}\subset T_{Hodge}^*(M_{t,n}(X_{\mathbb{Q}}))=: (W_{\cdot}\subset V_{\mathbb{Q}}, F^{\cdot}\subset V_{\mathbb{C}})
 \subset \Psi_t ^{Hodge}( H^n(X_{\mathbb{C}}/C_{\mathbb{C}}))\footnote{Here
 $T_{Hodge}^*(M_{t,n}(X_{\mathbb{Q}}))$ denotes the Hodge structure dual to
 $T^{Hodge}(M_{t,n}(X_{\mathbb{Q}}))$.} .$$
Hence, $W_0 L^{Hodge}  \otimes \mathbb{Q}  \simeq \Lambda ^* \otimes
\mathbb{Q}$ and  $W_{-1}M_{t,n}(X_{\mathbb{Q}}) =
W_{-2}M_{t,n}(X_{\mathbb{Q}})$. We claim that the image of the
embedding $W_2 L^{Hodge}/W_0 L^{Hodge} \otimes
\mathbb{Q}\hookrightarrow (T_* \otimes \mathbb{Q})^*$ is
$Gal(\overline{\mathbb{Q}}/\mathbb{Q})$-invariant. Indeed, this is
clear from the commutative diagram

   $$
\def\normalbaselines{\baselineskip20pt
\lineskip3pt  \lineskiplimit3pt}
\def\mapright#1{\smash{
\mathop{\to}\limits^{#1}}}
\def\mapdown#1{\Big\downarrow\rlap
{$\vcenter{\hbox{$\scriptstyle#1$}}$}}
\begin{matrix}
 L^{Hodge}_{\mathbb{Q}} \otimes \mathbb{Q}_p &    \stackrel{ }{\longrightarrow}   &   V_{\mathbb{Q}}\otimes \mathbb{Q}_p \cr
\mapdown{\simeq }  &  &\mapdown{\simeq}    \cr
  L^{et}\otimes _{\mathbb{Z}_p} \mathbb{Q}_p:= Im \, N^{n-1}_{et}  &    \stackrel{  }{\longrightarrow} &
T^{*et}_{\mathbb{Q}_p}(M_{t,n}(X_{\mathbb{Q}}))
\end{matrix}
 $$
since all the arrows in the bottom row are morphisms of
$Gal(\overline{\mathbb{Q}}/\mathbb{Q})$-modules. It follows that
there exists a unique quotient  $L^{mot} \in {\cal M}_1(\mathbb{Q};
\mathbb{Q})$, $\gamma: M_{t,n}(X_{\mathbb{Q}})\twoheadrightarrow
L^{mot} $ which fits into the following diagram
    $$
\def\normalbaselines{\baselineskip20pt
\lineskip3pt  \lineskiplimit3pt}
\def\mapright#1{\smash{
\mathop{\to}\limits^{#1}}}
\def\mapdown#1{\Big\downarrow\rlap
{$\vcenter{\hbox{$\scriptstyle#1$}}$}}
\begin{matrix}
 T_{Hodge}^*(L^{mot})&    \stackrel{ }{\longrightarrow}   &
  T_{Hodge}^*(M_{t,n}(X_{\mathbb{Q}})) \cr
\mapdown{\simeq }  &  &\mapdown{Id}    \cr
  L^{Hodge} &    \stackrel{  }{\longrightarrow} &
T_{Hodge}^*(M_{t,n}(X_{\mathbb{Q}} ))
\end{matrix}
 $$
Observe that the operator $N$ descends to $L^{mot}$ and
\begin{equation}\label{mon}
   N: (W_0 L^{mot}  /W_{-2} L^{mot})\otimes \mathbb{Q}(1)\stackrel{\simeq }{\longrightarrow}  W_{-2} L^{mot}\otimes \mathbb{Q} .
\end{equation}
Finally, we have from (\ref{c2})  a canonical isomorphism
$T^{*et}_{\mathbb{Q}_p}(L^{mot}) \simeq L^{et}\otimes
_{\mathbb{Z}_p} \mathbb{Q}_p$ of
 $Gal(\overline{\mathbb{Q}}/\mathbb{Q})$-modules.

Let $[L^{mot}, N^{-1}]\in Ext^1_{{\cal M}_1(\mathbb{Q};
\mathbb{Q})}(\mathbb{Q}(0), \mathbb{Q}(1))$ be the class of the
extension
\begin{equation}\label{motext}
 0 \to  W_{-2} L^{mot} \to L^{mot} \to W_0 L^{mot}  /W_{-2} L^{mot} \to 0
\end{equation}
 composed with $N^{-1}$ from (\ref{mon}) and let $\kappa$ be the corresponding (by (\ref{tateext})) element in  $\mathbb{Q}^* \otimes \mathbb{Q}$.
The functor  $T_{Hodge}^*$ takes this extension to the class
$[L^{Hodge}, N^{-1}_B] \in   Ext^1_{MHS}(\mathbb{Q}(0),
\mathbb{Q}(1)) \simeq \mathbb{C}^* \otimes \mathbb{Q}$. The latter
class is equal to $q'_{\mathbb{C}}(0)$. It follows that $
q'_{\mathbb{C}}(0)= \kappa ^{-1} $.

If we pull back the extension (\ref{motext}) on $spec \,
\mathbb{Q}_p$ and then apply the etale realization functor
   $T_{\mathbb{Q}_p}^*:  {\cal M}_1(\mathbb{Q}_p; \mathbb{Q}) \to Rep_{\mathbb{Q}_p}(\Gamma)$
   we get the extension
$[L^{et}\otimes \mathbb{Q}_p , N^{-1}_{et}]$ equivalent (by {\bf
4.2} ,\, 3)) to the one obtained from $[L^{FL}\otimes \mathbb{Q}_p,
N^{-1}_{DR}]$ by applying the Fontaine-Laffaille functor $U$. Hence
$ q'_{\mathbb{Z}_p}(0)= \kappa ^{-1} $, and we are done.

\begin{rem} The above argument shows that for any
family $X_{\mathbb{Q}} \to C_{\mathbb{Q}} $ over $\mathbb{Q}$ with a
maximal degeneracy point at $a\in \overline
C_{\mathbb{Q}}(\mathbb{Q})$
$$q_{\mathbb{C}}\in (\mathbb{Q}((t)))^*\otimes \mathbb{Q}.$$
\end{rem}

b) Assume that the Betti monodromy is small i.e.
\begin{equation}\label{monbetti}
N_B: W_2 L^{Hodge}_{\mathbb{Z}}/ W_0 L^{Hodge}_{\mathbb{Z}} \simeq
W_0 L^{Hodge}_{\mathbb{Z}}(-1)
\end{equation}
We have to show that, for any prime $p$ in $S$,
\begin{equation}\label{monrham}
N_{DR}: W_2 L^{FL}/ W_0 L^{FL} \to  W_0 L^{FL}(-1)
\end{equation}
is an isomorphism as well. Indeed, by   (by  {\bf 4.2} ,\, 3)) the
functor $U$ takes the morphism (\ref{monrham}) to
$$N_B\otimes Id:
(W_2 L^{Hodge}_{\mathbb{Z}}/ W_0 L^{Hodge}_{\mathbb{Z}})\otimes
\mathbb{Z}_p \simeq  W_0 L^{Hodge}_{\mathbb{Z}} \otimes
\mathbb{Z}_p(-1)  .$$ The claim follows.

c)  Let ${\cal E}$ be the quotient of $H^n_{log}(\overline
X_S/\overline C_S)$ modulo torsion, and let $\overline{\cal F}^n
\subset {\cal E}$, $ {\cal W}_0 \subset {\cal E}$ be the Hodge and
monodromy filtrations ({\bf 3.3}). As we explained in {\it loc.
cit.} the Poincare duality identifies the line bundle
$\overline{\cal F}^n$ with the dual to $ {\cal W}_0$. It is also
shown there that $ {\cal W}_0$ is generated by a parallel section
$e^0 \in {\cal W}_0^{\nabla}$. It suffices to prove the claim for a
single nonvanishing section $\omega \in \overline {\cal F}^n$. Let
us choose $\omega$ such that $(e^0, \omega)=0$. Then the integral
$$\frac{1}{(2\pi i)^n} \int
_{\delta _1} \omega $$ is a constant function on $D_{\mathbb{C}}$.
We have to show that the square of this constant is in
$\mathbb{Z}[N^{-1}]^*$. The following lemma does the job.
\begin{lemma}
Let $E$ (resp. $H$ ) be the torsion free part of
$a^*(H^n_{log}(\overline X_S/ \overline C_S))$ (resp.
$\Psi_t^{an}(R^n\pi_{\mathbb{C} *}^{an}\mathbb{Z}[N^{-1}])$ ) , and
let
 $$E \hookrightarrow E\otimes \mathbb{C}\simeq  H \otimes \mathbb{C} \hookleftarrow H$$
be the isomorphism from (\ref{derhamv}). Then the two
$\mathbb{Z}[N^{-1}]$-lattices
 $$(W_{0} E)^{\otimes ^2} \hookrightarrow (W_{0} E)^{\otimes ^2} \otimes \mathbb{C} \simeq  (W_{0} H)^{\otimes ^2}
\otimes \mathbb{C}\hookleftarrow (W_{0} H)^{\otimes ^2}$$ coincide.
\end{lemma}
\begin{proof}
Indeed, consider the monodromy paring
$$\Xi: <\cdot, \cdot >_{mon}: (W_{0} E)^{\otimes ^2} \otimes \mathbb{C}\to \mathbb{C},$$
$$  <x,y>_{mon}= <x, N^{-n}_{DR} y>_{DR} =\pm <x, N^{-n}_{B}y >_B,$$
where $N^{-n}_B= (-2\pi i)^{-n}N^{-n}_{DR}: W_{0} E \otimes
\mathbb{C} \to W_{2n} E \otimes \mathbb{C} $. The monodromy paring
takes $(W_{0} E)^{\otimes ^2} \otimes \mathbb{Q}$ and $ (W_{0}
H)^{\otimes ^2} \otimes \mathbb{Q}$ into $\mathbb{Q}\subset
\mathbb{C}$. Therefore, since $rk \, W_0 E =1$,   $(W_{0}
E)^{\otimes ^2} \otimes \mathbb{Q} = (W_{0} H)^{\otimes ^2} \otimes
\mathbb{Q}$. Moreover, to prove that $(W_{0} E )^{\otimes ^2}=(W_{0}
H )^{\otimes ^2}$, it is enough to show that $\Xi((W_{0} E
)^{\otimes ^2})=\Xi((W_{0} H )^{\otimes ^2})$. Since the pairings
$$<\cdot, \cdot>_{DR} : W_{0} E \otimes W_{2n} E \to \mathbb{Z}[N^{-1}], \, <\cdot, \cdot>_B : W_{0} H \otimes
W_{2n} H \to \mathbb{Z}[N^{-1}]$$ are perfect, the claim would
follow if we prove that, for any prime $p$ in $S$, the cokernels of
the maps $N_{DR}^{n}: W_{0}E \otimes \mathbb{Z}_p \to W_{2n} E
\otimes \mathbb{Z}_p, \, N_{et}^{n}: W_{0}H \otimes \mathbb{Z}_p \to
W_{2n} H \otimes \mathbb{Z}_p$ have the same order. This follows
from parts 1) and 3) in  {\bf 4.2}.
\end{proof}

\begin{center}
{\bf REFERENCES}
\end{center}

[A1] J. Ayoub, {\it Les six op\'erations de Grothendieck et le
formalisme des cycles \'evanescents dans le monde motivique}, 2006.
Available electronically at
http://www.institut.math.jussieu.fr/~ayoub/These/THESE.pdf

[A2] J. Ayoub, {\it Private communications.}

[BK] L. Barbieri-Viale, B. Kahn, {\it On the derived category of
1-motives, I}, arXiv:0706.1498v1 [math.AG].

[BK2] L. Barbieri-Viale, B. Kahn, {\it On the derived category of
1-motives, II}. In preparation.

[B] A. Beilinson, {\it How to glue perverse sheaves}, K-theory,
Arithmetic and Geometry, LNM 1289.

[BOV] A. Beilinson,  A. Otwinowska,  V. Vologodsky, {\it
      Motivic sheaves over a curve}, work in progress.

[COGP] P. Candelas, X. de la Ossa, P. Green, L. Parkes, {\it A pair
of Calabi-Yau manifolds as an exactly soluble superconformal field
theory}, Nuclear Physics, B359(1991) 21

[D1] P.Deligne, {\it Local behavior of Hodge structures at
infinity}, AMS/IP Studies in Advanced Mathematics, Volume 1, 1997.

[D2] P.Deligne, {\it Th\'eorie de Hodge 2}, Publ. Math. IHES 40
(1971).

[D3] P.Deligne, {\it Th\'eorie de Hodge 3}, Publ. Math. IHES 44
(1974).

[Fa] G. Faltings, {\it Crystalline cohomology and p-adic Galois
representations}, Algebraic Analysis, Geometry and Number Theory
(J.Igusa, ed.) (1989).

[Hu1] A. Huber, {\it Realization of Voevodsky's motives}, J.
Algebraic Geom. 9 (2000), no. 4.

[Hu2] A. Huber, {\it Corrigendum to: "Realization of Voevodsky's
motives" } J. Algebraic Geom. 13 (2004), no. 1.

[Il] L.Illusie, {\it Logarithmic spaces (according to K.Kato)},
Perspect. Math., 15.

[KSV] M. Kontsevich, A. Schwarz, V. Vologodsky, {\em  Integrality of
the instanton numbers  and p-adic B-model},  Physics Letters B 637
(2006).

[LY] B.Lian, S.T.-Yau, {\it Mirror Maps, Modular Relations and
Hypergeometric Series,1}, arXiv:hep-th/9507151 v1 27 July 1995.

[M] D.Morrison, {\it Mirror Symmetry and Rational Curves on Quintic
Threefolds: A Guide for Mathematicians}, J. Amer. Math. Soc. 6
(1993), no. 1.

[O] F. Orgogozo, {\it Isomotifs de dimension inf\'erieure ou \'egale
\`a un}, Manuscripta Math. 115 (2004), no. 3.

[Ste] J. Steenbrink, {\it Limits of Hodge structures}, Inv. Math. 31
(1976).

[Sti] J. Stienstra, {\it Ordinary Calabi-Yau-3 Crystals}, Calabi-Yau
varieties and mirror symmetry (Toronto, ON, 2001), Fields Inst.
Commun., 38, Amer. Math. Soc., Providence, RI, (2003).

[V] V. Voevodsky, {\it Triangulated category of motives over a
field}, in "Cycles, Transfers, and Motivic Homology Theories",
Annals of Mathematics Studies 143 (2000).

[Vol] V. Vologodsky, {\it The Albanese functor commutes with the
Hodge realization.}
\end{document}